\documentclass[a4paper,12pt]{amsart}
\usepackage{amssymb,amscd,amsfonts}
\usepackage{calrsfs}

\theoremstyle{plain}
\newtheorem{theo}{Theorem}[section]
\newtheorem{prop}[theo]{Proposition}
\newtheorem{lemm}[theo]{Lemma}
\newtheorem{coro}[theo]{Corollary}

\theoremstyle{definition}

\newtheorem{exam}{Example}

\theoremstyle{remark}
\newtheorem{rema}[theo]{Remark}
\newtheorem*{note}{Note}
\newtheorem*{claim}{Claim}

\renewcommand{\labelenumi}{(\roman{enumi})}

\renewcommand{\u}{u_i^I}

\newcommand{\V}{\mathcal V}

\newcommand{\field}[1]{\mathbb{#1}}
\newcommand{\C}{\field{C}}

\newcommand{\Q}{\field{Q}}
\newcommand{\R}{\field{R}}
\newcommand{\Z}{\field{Z}}

 \DeclareMathOperator{\Hom}{Hom}

\newcommand{\var}{\varphi}
\newcommand{\varv}{\varphi^v}

\newcommand{\hatvar}{\hat{\varphi}}
\newcommand{\hatvarv}{\hat{\varphi}^v}

\newcommand{\brvarv}{\breve{\varphi}^v}
\newcommand{\brvar}{\breve{\varphi}}

\newcommand{\brT}{\breve{T}}

\newcommand{\hatH}{\hat{H}}

\newcommand{\hatT}{\hat{T}}

\newcommand{\barh}{\bar{h}}

\newcommand{\sigmn}{\Sigma^{(n)}}
\newcommand{\sigmk}{\Sigma^{(k)}}
\newcommand{\sigmone}{\Sigma^{(1)}}

\def\l{\langle}
\def\r{\rangle}

\newcommand{\Proj}{\field{P}}

\newcommand{\img}{\sqrt{-1}}

\title{Elliptic genera, torus orbifolds and multi-fans; II}
\author{Akio Hattori}
\address{Graduate School of Mathematical Science, University of Tokyo,
Tokyo, Japan}
\email{hattori@ms.u-tokyo.ac.jp}
 
\begin{document}
\maketitle

\section{Introduction}
\label{sec:intro}
This article is a continuation of \cite{HM2}. 

Elliptic genera for manifolds introduced by Ochanine and other people 
has a remarkable feature called rigidity.
If the circle group acts non-trivially on a closed almost complex 
(or more generally stably almost complex) manifold whose first Chern 
class is divisible
by a positive integer $N$ greater than $1$, then its equivariant
elliptic genus of level $N$ is rigid, that is, it is a constant 
character of the circle group. It was conjectured by Witten \cite{W}
and proved by Taubes \cite{T}, Bott-Taubes \cite{BT} and Hirzebruch 
\cite{Hir}. Liu \cite{L}
found a simple proof using the modular property of elliptic genera. 

For stably almost complex orbifolds an invariant called 
orbifold elliptic genus is defined. 
Though the rigidity property does not hold in this case 
in the same form as the case of manifolds, versions of rigidity 
theorem can be formulated and shown to hold \cite{Hat}. The 
orbifold elliptic genus is also defined for multi-fans and 
is likewise expected to have nice properties \cite{HM2}. 
The notion of multi-fan is a generalization of that of fan in 
the theory of toric varieties. It was first introduced in \cite{M} 
and was further developed in \cite{HM1}. In \cite{HM2} 
we proved a rigidity theorem concerning the elliptic genus 
for complete simplicial multi-fans 
and gave an application to toric varieties 
to the effect that a non-singular complete toric variety 
of dimension $n$ whose 
canonical line bundle is linearly equivalent to a $T$-Cartier 
divisor of the form $nD$ is isomorphic to a certain 
projective space bundle over a projective line. 

In \cite{Hat} we proved 
rigidity theorems and vanishing theorems concerening orbifold elliptic genus 
for genenral almost complex orbifolds. In this paper 
we shall give analogues of those vanishing theorems for multi-fans 
(Theorem \ref{theo:brrigid}, Theorem \ref{theo:hatrigid} and 
theorem \ref{theo:varrigid}). 
As a special feature about multi-fans, vanishing of orbifold 
elliptic genus and modified orbifold elliptic genus holds 
under suitable assumptions. 

As an application it will be shown that a $\Q$-factorial 
complete toric variety of dimenson $n$ whose canonical divisor 
is linearly equivalent to a $T$-Cartier divisor of the form $nD$ 
is either non-singular or isomorphic to a certain weighted 
projective space (Proposition \ref{prop:N=n} and Corollary \ref{coro:N=n}). 
This result was already proved by O. Fujino \cite{Fuj} for 
projective toric varieties. The author is grateful to him 
for having communicated his result to the author. 

The paper is organized as follows. In Section 2 we recall some basic facts
about multi-fans from \cite{HM1} and \cite{HM2}. 
In Section 3 we define the elliptic genus, orbifold elliptic 
genus and modified orbifold elliptic genus of a pair of a 
multi-fan and a set of generating vectors and formulate main theorems. 
Properties of the $T_y$-genus, orbifold $T_y$-genus and modified orbifold 
$T_y$-genus of multi-fans are discussed. 
In Section 4 the proofs of the main results are given. 
Section 5 is devoted to applications to special types of 
multi-fans and, in particular, fans associated to toric 
varieties. 

The author wishes to thank M. Masuda, coauthor of the joint papers 
\cite{HM1} and \cite{HM2}. The collaboration with him 
was much profitable to the author for this work as well. He declined to 
be coauthor of this paper which in fact grew up from our collaboration.

\section{Multi-fans}\label{sec:multi}
The present paper depends heavily on its first part \cite{HM2}. 
In this section we shall summarize materials we need in the 
sequel from \cite{HM1} and \cite{HM2}. 

Let $L$ be a lattice of rank n.  
An $n$-dimensional 
\emph{simplicial multi-fan} in $L$ is a triple 
$\Delta=(\Sigma,C,w^{\pm})$. We shall call it simply 
a multi-fan in this paper. Here $\Sigma$ is an
augmented simplicial set, that is, 
$\Sigma$ is a simplicial set with empty set 
$*=\emptyset$ added as $(-1)$-dimenional simplex.
$\sigmk$ denotes the $k-1$ skeleton of $\Sigma$ so that 
$*\in \Sigma^{(0)}$. 
We assume that $\Sigma=\sum_{k=0}^n\sigmk$, and $\sigmn \not=\emptyset$.
$C$ is a map from $\sigmk$ into the set of $k$-dimensional 
strongly convex rational
polyhedral cones in the vector space $L_\R=L\otimes \R$ for 
each $k$ such that, if $J$ is a face of $I$, then $C(J)$ is
a face of $C(I)$. 
$w^{\pm}$ are maps $\sigmn \to \Z_{\ge 0}$. 
We set $w(I)=w^+(I)-w^-(I)$.
A vector $v\in L_\R$ will be called \emph{generic} if $v$ does
not lie on any linear subspace spanned by a cone in 
$C(\Sigma)$ of dimesnsion less than $n$. For a generic
vector $v$ we set $d_v=\sum_{v\in C(I)}w(I)$, where
the sum is understood to be zero if there is no such $I$.
We call a multi-fan $\Delta=(\Sigma,C,w^\pm)$ of dimension $n$ 
{\it pre-complete} if 
the integer $d_v$ is independent of
the choice of generic vectors $v$. We call this integer
the {\it degree} of $\Delta$ and denote it by $\deg(\Delta)$.

For each $K\in \Sigma$ we set
\[ \Sigma_K=\{ J\in\Sigma\mid K\subset J\}.\]
It inherits the partial ordering from $\Sigma$ and becomes an augmented
simplicial set where $K$ is the unique 
minimum element in $\Sigma_K$. 
Let $(L_K)_\R$ be the linear subspace of $L_\R$ generated by $C(K)$.
Let $L^K_\R$ be the quotient space of $L_\R$ by $(L_K)_\R$ and 
$L^K$ the image of $L$ in $L^K_\R$. $L^K_\R$ is identified with
$L^K\otimes \R$. 
For $J\in \Sigma_K$ we define $C_K(J)$ to be
the cone $C(J)$ projected on $L^K_\R$. 
We define two functions 
\[ {w_K}^\pm\colon  \Sigma_K^{(n-|K|)}\subset \Sigma^{(n)} \to \Z_{\ge 0}\]
to be the restrictions of $w^\pm$ to $\Sigma_K^{(n-|K|)}$. The triple 
$\Delta_K:=(\Sigma_K,C_K,{w_K}^\pm)$ is a multi-fan in $L^K$ 
and is called the  
\emph{projected multi-fan} with respect to $K\in \Sigma$. 
If $K=\emptyset$ then $\Delta_K =\Delta$. 
A pre-complete multi-fan $\Delta=(\Sigma,C,w^\pm)$ is 
said to be 
\emph{complete} if the projected multi-fan $\Delta_K$ is pre-complete for any 
$K\in \Sigma$. A multi-fan is complete if and only if the projected
multi-fan $\Delta_J$ is pre-complete for any $J\in \Sigma^{(n-1)}$.

Let $M$ be an oriented closed manifold of dimension $2n$ with an 
effective action of an $n$-dimensional torus $T$. We assume further 
that the fixed point set $M^T$ is not empty. There is a 
finite number of subcircles of $T$ such that the fixed point
set of each subcircle has codimension 2 components. Let $\{M_i\}_{i=1}^r$
be those components which have non-empty intersection with $M^T$.
We call $M$ a \emph{torus manifold} if a preferred orientation of 
each $M_i$ is given. The $M_i$ are called \emph{characteristic submanifolds}. 
A complete multi-fan 
$\Delta(M)=(\Sigma(M),C(M),w^\pm(M))$
in the lattice $H_2(BT)$ is associated with $M$ in a canonical way, where 
$BT$ is the classifying space of $T$. 
If $K\in \Sigma(M)$ then the projected multi-fan $\Delta(M)_K$ is closely 
related to the multi-fan associated with $M_K=\cap_{i\in K}M_i$, where
$M_K$ is regarded as a union of torus manifolds. 

Let $\Delta=(\Sigma,C,w^\pm)$ be a multi-fan in $L$. 
If $T$ denotes the torus $L_\R/L$, then $L$ can be canonically
identified with $H_2(BT)$.
Then there is a unique
primitive vector $v_i\in L=H_2(BT)$ which generates the cone $C(i)$
for each $i\in \sigmone$. $\Delta$ is called \emph{non-singular} if
$\{v_i\mid i\in I\}$ is a basis of the lattice $L=H_2(BT)$ for each 
$I\in \sigmn$. The multi-fan associated with a torus manifold
is a complete non-singular multi-fan.

It is sometimes more convenient to consider a set of vectors
${\mathcal V}=\{v_i\in L\}_{i\in \sigmone}$ such that each $v_i$ generates
the cone $C(i)$ in $L_\R$ but is not necessarily primitive. 
This is the case for multi-fans associated with 
torus orbifolds. A torus orbifold $X$ is a closed
oriented orbifold of even dimension with an effective action of 
a torus of half the dimension of the orbifold $X$ with some additional
condition. 
A set of codimension 2 suborbifolds $X_i$ called \emph{characteristic 
suborbifolds} is similarly defined
as in the case of torus manifolds. To each subcircle $S_i$
which fixes $X_i$ pointwise there is a unique finite covering 
$\tilde{S}_i$ and an effective action of $\tilde{S}_i$ on the 
orbifold cover of each fiber of the normal bundle.
This defines a vector $v_i$ in $\Hom(S^1,T)=H_2(BT)=L$ as before.
In this way a multi-fan $\Delta(X)$ and a set of vectors
${\mathcal V}(X)=\{v_i\}_{i\in \sigmone}$ are associated to
the torus orbifold $X$. 

Hereafter multi-fans are assumed to be complete 
and we shall always consider the pair of a multi-fan $\Delta$ 
and a set of vectors
${\mathcal V}=\{v_i\in L\}_{i\in \sigmone}$ as above. 
In case $\Delta$ is non-singular it is further assumed 
that all the $v_i$ are primitive.
If $I$ is in $\sigmn$, then $\{v_i\}_{i\in I}$ becomes
a basis of vector space $L_\R$.
In case $\Delta$ is non-singular it is a basis of the lattice
$L$. In general, for $I\in \sigmn$, we define $L_{I,\V}$ to be 
the sublattice of $L$
generated by $\{v_i\}_{i\in I}$. 

Let $L_{I,\V}^*$ be the dual lattice of
$L_{I,\V}$ and $\{u_i^I\}$ the basis of $L_{I,\V}^*$ dual 
to $\{v_i\}_{i\in I}$. 
We identify $L_{I,\V}^*$ with the lattice in $L_\R^*$
given by 
\[ \{u\in L_\R^*\mid \langle u,v\rangle \in \Z,
   \ \ \text{for any}\ \  v\in L_{I,\V}\},\] 
where $\l u,v\r$ is the dual pairing. 
For $h\in L/L_{I,\V}$ and $u\in L_{I,\V}^*$ we define 
\[ \chi_I(u,h)=e^{2\pi\sqrt{-1}\langle u,v(h)\rangle},\]
where $v(h)\in L$ is a representative of $h$. If one fixes $u$,
$h \mapsto \chi_I(u,h)$ gives a character of the group $L/L_{I,\V}$. 

The dual lattice
$L^*=H^2(BT)\subset H^2(BT;\R)$ is canonically identified 
with $\Hom(T,S^1)$. The latter is embedded in the character ring
$R(T)$. In fact $R(T)$ can be considered as the group ring $\Z[L^*]$ of 
the group $L^*=\Hom(T,S^1)$. It is convenient to write the element 
in $R(T)$ corresponding to $u\in H^2(BT)$ by $t^u$. The 
homomorphism $v^*:R(T)\to R(S^1)=\Z[t,t^{-1}]$ induced by an element 
$v\in H_2(BT)=\Hom(S^1,T)$ can be written in the form
\[ v^*(t^u)=t^{\langle u,v\rangle}, \]
where $t^m\in R(S^1)$ is such that $t^m(g)=g^m$ for $g\in S^1$.

More generally, set $L_\V=\bigcap_{I\in \sigmn}L_{I,\V}$, and 
let $L_{\V}^*$ be the dual lattice of $L_\V$. 
$L_{\V}^*$ contains all $L_{I,\V}^*$ and is generated by all the 
$u_i^I$'s. The group ring $\Z[L_{\V}^*]$ 
contains $\Z[L^*]=R(T)$ and has a basis $\{t^u|u\in L_{\V}^*\}$ with 
multiplication determined by the addition in $L_{\V}^*$:
\[  t^ut^{u'}=t^{u+u'}. \]
If $v$ is a vector in $L_\V$, then $v$ determines a homomorphism
$v^*:\Z[L_{\V}^*]\to R(S^1)=\Z[t,t^{-1}]$ sending $t^u$ to
$t^{\langle u,v\rangle}$. If we vary $v$ then 
$v^*(t^u)$ determines $t^u$. 

 Similarly if $v_1$ and $v_2$ are vectors 
in $L$ they define a homomorphism from a $2$-dimensional 
torus $T^2$ into $T$ and induce a homomorphism 
$(v_1,v_2)^*:\Z[L^*]\to R(T^2)=\Z[t_1,t_1^{-1},t_2,t_2^{-1}]$ 
defined by 
\[ (v_1,v_2)^*(t^u)=t_1^{\l u,v_1\r}t_2^{\l u,v_2\r}. \]
If $v_1$ and $v_2$ belong to $L_\V$, then $(v_1,v_2)^*$ extends to 
a homomorphism $\Z[L_\V^*]\to R(T^2)$. 
We define the equivariant cohomology $H_T^*(\Delta)$ of a complete 
multi-fan $\Delta$ as the face ring of the simplicial complex $\Sigma$.
Namely let $\{x_i\}$ be indeterminates indexed by $\sigmone$, and 
let $R$ be the polynomial ring over the integers generated 
by $\{x_i\}$. We denote by $\mathcal{I}$ the ideal in $R$
generated by monomials
$\prod_{i\in J} x_i$ such 
that $J\notin \Sigma$. $H_T^*(\Delta)$ is by definition
the quotient $R/\mathcal{I}$. 
We regard $H^2(BT)$ as a submodule of $H_T^2(\Delta)$ by the formula
\begin{equation}\label{eq:structure}
 u=\sum_{i\in \sigmone}\langle u,v_i\rangle x_i .
\end{equation} 
This determines an $H^*(BT)$-module structure of $H_T^*(\Delta)$.
It should be noticed that this module structure depends on
the choice of ${\mathcal V}$ as above.

For each $I\in\sigmn$ we define the restriction homomorphism
$\iota_I^* :H_T^2(\Delta)\to L_{\V}^*$ by
\begin{equation*}\label{eq:restrict}
 \iota_I^*(x_i)=\begin{cases}
            u_i^I & \text{for}\ i\in I\\
            0     & \text{for}\ i\notin I.
            \end{cases}
\end{equation*}            
$\iota_I^*| H^2(BT)$ is 
the identity map for any $I$, and $\sum_{I\in \sigmn}\iota_I^*$ is 
injective. Note that, if $\Delta$ is non-singular, then $\iota_I^*$ maps
$H_T^2(\Delta)$ into $H^2(BT)$.

\begin{lemm}\label{lemm:multiplicity bis}
For any $x=\sum_{i\in \sigmone}c_ix_i\in H_T^2(\Delta),\ c_i\in \Z$,
the element 
\[ \sum_{I\in\sigmn}\frac{w(I)}{|L/L_{I,\V}|}\sum_{h\in L/L_{I,\V}}
           \frac{\chi_I(\iota_I^*(x),h)t^{\iota_I^*(x)}}
         {\prod_{i\in I}(1-\chi_I(u_i^I,h)^{-1}t^{-u_i^I})}\]
in $\C[L_{\V}^*]$ actually belongs to $R(T)$.
\end{lemm}
See Lemma 2.1 in \cite{HM2}. 

We also use an extended version of this Lemma. 
Let $K\in \sigmk$ and let $\Delta_K=(\Sigma_K,C_K,w_K^{\pm})$ be 
the projected multi-fan.
If $I\in\Sigma^{(l)}$ contains $K$, then $I$ is considered 
as lying in $\Sigma_K^{(l-k)}$. In order to avoid some notational confusions
we introduce the link $\Sigma'_K$ of $K$ in $\Sigma$.
It is a simplicial set consisiting simplices $J$ such that 
$K\cup J\in \Sigma$ and $K\cap J=\emptyset$. There is an isomorphism 
from $\Sigma'_K$ to $\Sigma_K$ sending $J\in {\Sigma'_K}^{(l)}$ to 
$K\cup J\in \Sigma_K^{(l)}$. Let $K\ast\Sigma'_K$ be the join of 
$K$ (regarded as a simplicial set) and $\Sigma'_K$. Its simplices are
of the form $J_1\cup J_2$ with $J_1\subset K$ and $J_2\in \Sigma'_K$. 
The torus $T^K$ corresponding to $\Delta_K$ 
is a quotient of $T$. We consider the polynomial 
ring $R_K$ generated by $\{x_i\mid i\in K\cup{\Sigma'_K}^{(1)}\}$ 
and the ideal $\mathcal{I}_K$ generated by monomials
$\prod_{i\in J} x_i$ such 
that $J\notin K\ast\Sigma'_K$. We define the equivariant cohomology 
$H_T^*(\Delta_K)$ of
$\Delta_K$ with respect to the torus $T$ as the quotient ring
$R_K/\mathcal{I}_K$. Note that $H_T^*(\Delta_K)$ is defferent from 
$H_{T^K}^*(\Delta_K)$. 

$H^2(BT)$ is regarded as a submodule of
$H_T^2(\Delta_K)$ by a formula similar to (\ref{eq:structure}).
This defines
an $H^*(BT)$-module structure on $H_T^*(\Delta_K)$ .
The projection $H_T^2(\Delta)\to H_T^2(\Delta_K)$ 
is defined by sending $x_i$ to $x_i$ for $i\in K\cup{\Sigma'_K}^{(1)}$ and
putting $x_i=0$ for $i\notin K\cup{\Sigma'_K}^{(1)}$.  
The restriction homomorphism
$\iota_I^*:H_T^2(\Delta_K) \to L_{\V}^*$ is also defined for 
$I\in \Sigma_K^{(n-k)}$ by $\iota_I^*(x_i)=u_i^I$. If $X$ is a torus 
orbifold, then 
$H_T^*(\Sigma(X)_K)$ is related to the equivariant cohomology
$H_T^*(X_K)$ with respect to the group $T$ (not with respect to $T_K$).

\begin{theo}\label{theo:multiplicity}
Let $\Delta$ be a complete simplicial multi-fan.
Let $x=\sum_{i\in K\cup{\Sigma'_K}^{(1)}}c_ix_i\in H_T^2(\Delta_K)$ be as 
above with all $c_i$ integers. Then the expression 
\begin{equation*}
 \sum_{I\in \Sigma_K^{(n-k)}}\frac{w(I)}{|L/L_{I,\V}|}\sum_{h\in L/L_{I,\V}}
 \frac{\chi_I(\iota_I^*(x),h)t^{\iota_I^*(x)}}
 {\prod_{i\in I\setminus K}(1-\chi_I(u_i^I,h)^{-1}t^{-u_i^I})}.
\end{equation*}
belongs to $R(T)$.
\end{theo}
See Theorem 2.3 in \cite{HM2}. Applying $(v_1,v_2)^*$ to the 
above expression we get

\begin{coro}\label{coro:multiplicity}
Let $v_1$ and $v_2$ be generic vectors in $L$ such that 
$\l \u,v_1\r$ and $\l\u,v_2\r$ are intergers for all $I\in \Sigma_K^{(n-k)}$, 
and let $x=\sum_{i\in K\cup{\Sigma'_K}^{(1)}}c_ix_i\in H_T^2(\Delta_K)$ 
with all $c_i$ integers. Then 
\begin{equation*}
\sum_{I\in \Sigma_K^{(n-k)}}\frac{w(I)}{|L/L_{I,\V}|}\sum_{h\in L/L_{I,\V}}
 \frac{\chi_I(\iota_I^*(x),h)
  t_1^{\l\iota_I^*(x),v_1\r}t_2^{\l\iota_I^*(x),v_2\r}}
   {\prod_{i\in I\setminus K}
  (1-\chi_I(u_i^I,h)^{-1}t_1^{\l-u_i^I,v_1\r}t_2^{\l-u_i^I,v_2\r})} 
\end{equation*}
belongs to $R(T^2)=\Z[t_1,t_1^{-1},t_2,t_2^{-1}]$.
\end{coro}
See Corollary 2.4 of \cite{HM2}. Note that if $v_1$ and 
$v_2$ belong to $L_\V$ then they satisfy the condition in 
Corollary \ref{coro:multiplicity}. 

Let $\Delta$ be a complete multi-fan 
in a lattice $L$ and $\V=\{v_i\}_{i\in \sigmone}$ a set of prescibed 
edge vectors as before. An $H^*(BT)$-module structure of 
$H_T^*(\Delta)$ is defined by \eqref{eq:structure}. The class 
\[ \sum_{i\in \sigmone}x_i \in H_T^2(\Delta) \]
will be called the \emph{equivariant first Chern class} of the 
pair $(\Delta,\V)$, and 
will be denoted by $c_1^T(\Delta,\V)$. When $\Delta$ is non-singular,
$\V$ consists of primitive vectors which is determined by $\Delta$ by
our convention. In this case we simply write $c_1^T(\Delta)$ 
and call it the \emph{equivariant first Chern class} of the 
non-singular multi-fan $\Delta$.

The image of $c_1^T(\Delta,\V)$ in
$H^2(\Delta,\V):=H_T^2(\Delta)/H^2(BT)$ is called the first Chern class
of $(\Delta,\V)$ and is denoted by $c_1(\Delta,\V)$.
Let $N>1$ be an integer. The first Chern class $c_1(\Delta,\V)$ is
divisible by $N$ in $H^2(\Delta,\V)$ if and only 
if $c_1^T(\Delta,\V)$ is of the form
\begin{equation}\label{eq:Ndivide}
 c_1^T(\Delta,\V)=Nx +u,\ x\in H_T^2(\Delta),\ u\in H^2(BT). 
\end{equation} 
We set $u^I=\iota_I^*(c_1^T(\Delta,\V))=\sum_{i\in I}u_i^I\in L_{I,\V}^*$.
Note that $u^I$ does not belong to $L^*=H^2(BT)$ in general.

\begin{lemm}\label{lemm:Ndivisible}
The following three conditions are equivalent:
\begin{enumerate}
\item the first Chern class $c_1(\Delta,\V)$ is
divisible by $N$, 
\item $u^I \bmod N $ regarded as an element of 
$L_{\V}^*/NL_{\V}^*$ is independent of $I\in \sigmn$ and 
belongs to the image of $L^*=H^2(BT)$,
\item there is an element $u\in H^2(BT)$ such that 
$\langle u,v_i\rangle =1\bmod N$ for all $i\in \sigmone$.
\end{enumerate}
\end{lemm}
See Lemma 4.1 in \cite{HM2}.

\begin{rema}\label{rem:face ring}
Let $M$ be a torus manifold
and $\Delta(M)$ its associated multi-fan. 
Put $\hat{H}_T^2(M)=H_T^2(M)/S\text{-torsion}$, where 
$S$ is the subset of $H^*(BT)$ multiplicatively generated by 
non-zero elements in $H^2(BT)$. In \cite{M} 
it was shown that there is a canonical
embedding of $H_T^2(\Delta(M))$ in $\hat{H}_T^2(M)$, and, in case
$M$ is a stably almost complex torus manifold, 
$c_1^T(M)\in H_T^2(M)$ descends to $c_1^T(\Delta(M))\in H_T^2(\Delta(M))$. 
It follows that,
if $M$ is a stably almost complex torus manifold and $c_1(M)$ is divisible by
$N$, then $c_1(\Delta(M))$ is also divisible by $N$. Even if $X$ is
a stably almost complex orbifold it can be shown that 
$c_1^T(X)\in H_T^2(X;\R)$ descends to 
$c_1^T(\Delta(X),\V(X))\in H_T^2(\Delta(X))\otimes\R$. But the divisibility of
the first Chern class has no meaning with real coefficients. 
We have to work with orbifold cohomology theory with integer coefficents, 
cf. Remark 3.2 in \cite{Hat}.
\end{rema} 

An element $x=\sum_{i\in \sigmone}c_ix_i$ is called \emph{$T$-Cartier} 
if $\iota_I^*(x)=\sum_{i\in \sigmone}c_i\u$ belongs to $L^*$ for 
each $I\in \sigmn$. Note 
that this definition depends on the choice of $\V$. The 
equivariant first Chern class $c_1^T(\Delta,\V)$ is 
$T$-Cartier if and only if $u^I=\sum_{i\in I}u_i^I$ belongs 
to $L^*$ for each $I\in \sigmn$. The Chern class $c_1(\Delta,\V)$ is said 
to be \emph{$T$-Cartier divisible} by $N$ if $x$ is $T$-Cartier in 
\eqref{eq:Ndivide}. In this case $c_1^T(\Delta,\V)$ is also 
$T$-Cartier. 
\begin{rema}
Let $X$ be a complete $\Q$-factorial toric variety. Its fan can be 
considered as a simplicial multi-fan  
$\Delta(X)=(\Sigma(X),C(X),w^\pm(X))$ in the lattice 
$L=\Hom(S^1,T)=H_2(BT)$ with $w^+(X)=1,w^-(X)=0$ 
for any $I\in \sigmn(X)$. We take the primitve vector $v_i$ of the 
$1$-dimensional cone $C(i)$ for each $i\in\sigmone$ and set $\V=\{v_i\}$. 
Let $D_i$ be the divisor corresponding to the cone $C(i)$. In our language 
it is the divisor corresponding to the characteristic suborbifold $X_i$. 
Every $T$-Weil divisor $D$ is written uniquely in the form
\[ D=\sum_{i\in \sigmone}c_iD_i,\ c_i\in \Z, \]
and $D$ is $T$-Cartier if and only if $\sum_{i\in \sigmone}c_i\u\in L^*$ for 
each $I\in \sigmn$, cf. \cite{Ful}. This suggests the above definition 
of $T$-Cartier elements in $H_T^*(\Delta)$. 
Each $u\in L^*$ determines a rational function on $X$. The corresponding 
Cartier divisor is denoted by $[t^u]$. 
Let $Div_T(X)$ be the group of $T$-Cartier 
divisors, $Pic(X)$ the group of line bundles on $X$ and $A_{n-1}(X)$ 
the group of all Weil divisors. Then 
there is a morphism of exact sequnces  
\[ \begin{CD}
 0 @>>> L^* @>>> Div_T(X) @>>> Pic(X) @>>> 0 \\
 @.     @|        @VVV         @VVV     \\
 0 @>>> L^* @>>> \bigoplus_{i\in \sigmone}\Z\cdot D_i @>>> A_{n-1}(X) 
 @>>> 0 \end{CD} \]
where the map $L^*=\Hom(T,S^1)\to Div_T(X)$ is given by 
\begin{equation}\label{eq:principal}
 u\mapsto [div(t^u)]=\sum_i\l u,v_i\r D_i.
\end{equation} 
See \cite{Ful}. Comparing \eqref{eq:structure} and \eqref{eq:principal} 
we see that the group of all $T$-Weil divisors 
$\bigoplus_{i\in \sigmone}\Z\cdot D_i$ can be identified with 
$H_T^2(\Delta)$ and $T$-Cartier divisors with $T$-Cartier elements. 
\end{rema}

\section{Elliptic genera of multi-fans}
\label{sec:elliptic of multi}
Let $\Delta$ be a complete simplicial multi-fan in a lattice $L$ 
and $\V=\{v_i\}_{i\in\sigmone}$
a set of prescribed edge vectors as in Section 2. We first recall 
the definition of 
(equivariant) elliptic genus $\var(\Delta,\V)$ and 
(equivariant) orbifold elliptic genus $\hatvar(\Delta,\V)$ of 
the pair $(\Delta,\V)$ after \cite{HM2}. They are defined in 
such a way that, for an almost complex torus orbifold $X$, 
$\var(\Delta(X),\V(X))$ and $\hatvar(\Delta(X),\V(X))$ coincide 
with the elliptic genus $\var(X)$ and the orbifold elliptic 
genus $\hatvar(X)$ of $X$ respectively. 

We consider the function $\phi(z,\tau,\sigma)$ of 
$z,\sigma$ in $\C$ and $\tau$ in the upper half plane $\mathcal{H}$ 
given by the following formula 
\[ \phi(z,\tau,\sigma)
=\zeta^{-\frac12}\frac{1-\zeta t}{1-t}
\prod_{k=1}^{\infty}\frac{(1-\zeta tq^n)(1-\zeta^{-1}t^{-1}q^n)}
 {(1-tq^n)(1-t^{-1}q^n)}, \]
where $t=e^{2\pi\img z},q=e^{2\pi\img\tau}$ and $\zeta =e^{2\pi\img\sigma}$. 
Note that $|q|<1$ and $\phi$ is a meromorphic function of $z,\tau$ 
and $\sigma$. 
We let act the group $SL_2(\Z)$ on $\C\times \mathcal{H}$ by 
\[ A(z,\tau)=(\frac{z}{c\tau+d}, A\tau)=
 (\frac{z}{c\tau+d},\frac{a\tau+b}{c\tau+d}),\ \ 
 A=\begin{pmatrix} a & b\\ c & d \end{pmatrix}\in SL_2(\Z). \] 
Then 
$\phi$ satisfies the following transformation formulae, cf. \cite{HM2}:
\begin{equation}\label{eq:transformula1}
 \begin{split}
 \phi(A(z,\tau),\sigma)&=
  e^{\pi\img c(2z\sigma +(c\tau +d)\sigma^2)}\phi(z,\tau,(c\tau +d)\sigma), 
  \\
  \phi(z+m\tau +n,\tau,\sigma)&=e^{-2\pi\img m\sigma}\phi(z,\tau,\sigma)=
  \zeta^{-m}\phi(z,\tau,\sigma).
 \end{split}
\end{equation}

In the sequel we fix the set $\V$ and put $H_I=L/L_{I,\V}$. 
Let $v\in L_{\V}$ be a generic vector. 
The (equivariant) \emph{elliptic genus} $\varv(\Delta,\V)$ 
along $v$ and 
the (equivariant) \emph{orbifold elliptic genus} 
$\hatvarv(\Delta,\V)$ along $v$
of the pair $(\Delta,\V)$ are defined by
\begin{equation*}\label{eq:varv}
 \varv(\Delta,\V)=\sum_{I\in\sigmn}\frac{w(I)}{|H_I|}
  \sum_{h\in H_I}
   \prod_{i\in I}\phi(-\langle u_i^I,zv+v(h)\rangle,\tau,\sigma),
\end{equation*} 
and
\begin{equation}\label{eq:hatvarv}
 \begin{split}
 \hatvarv(\Delta,\V) 
 = &\sum_{I\in\sigmn}
   \frac{w(I)}{|H_I|}\cdot \\
   &\sum_{(h_1,h_2)\in H_I\times H_I}
  \prod_{i\in I}\zeta^{\l u_i^I, v(h_1)\r}
 \phi(-\l u_i^I, zv-\tau v(h_1)+v(h_2)\r, \tau, \sigma).
\end{split} 
\end{equation}
where $v(h),v(h_1),v(h_2)\in L$ are representatives of $h,h_1,h_2\in H_I$
respectively.  
The above expressions give well-defined functions independent of 
the choice of representatives
$v(h),v(h_1),v(h_2)$ as is easily seen from \eqref{eq:transformula1}.
They are meromorphic functions in the variables $z,\tau,\sigma$ 
and sometimes written as $\varv(\Delta,\V;z,\tau,\sigma)$
and $\hatvarv(\Delta,\V;z,\tau,\sigma)$ to 
emphasize the variables. 

For each $K\in \sigmk$ with $k>0$ let $L_K$ be the kernel 
of the projection map
$L\to L^K$ and let $L_{K,\V}$ be the sublattice of $L_K$
generated by $v_i \in K$. We set 
 $H_K=L_K/L_{K,\V}$. If $J\subset K$ then we have
$ L_J\cap L_{K,\V}=L_{J,\V} $, and hence 
$H_J$ is canonically 
embedded in $H_K$. We set
\[\hat{H}_K=H_K\setminus \bigcup_{J\subsetneqq K}H_J .\]
The subset $\hat{H}_K$ is characterized by
\begin{equation}\label{eqn:0}
 \hat{H}_K=\{h\in H_K\mid \langle u_i^K,v(h)\rangle \not\in \Z 
\quad \text{for any $i\in K$} \}, 
\end{equation}
where $\{u_i^K\}$ is the basis of $L_{K,\V}^*$
dual to the basis $\{v_i\}_{i\in K}$ of $L_{K,\V}$ and $v(h)\in L_K$ is a 
representative of $h\in H_K$. 
For the minimum element $*=\emptyset\in \Sigma^{(0)}$
we set $\hat{H}_*=H_*=0$. 

If $K$ is contained in $I\in \sigmn$, then the canonical map
$L_{I,\V}^*\to L_{K,\V}^*$ sends $u_i^I$ to $u_i^K$ for $i\in K$
and to $0$ for $i\in I\setminus K$. Therefore, if $h$ is in $H_K$, 
then $\l u_i^I,v(h)\r=0$ for $i\in I\setminus K$, and 
$\l u_i^I,v(h)\r=\l u_i^K,v(h)\r$ for $i\in K$. Here $v(h)\in L_K$ is 
regarded as lying in $L$. Then 
$\hatvarv(\Delta,\V)$ can also be written in the following form 
which is sometimes useful.

\begin{equation}\label{eq:hatvarv2}
 \begin{split}
 \hatvarv(\Delta,\V) 
 = \sum_{k=0}^n&\sum_{K\in\sigmk,h_1\in \hat{H}_K}
   \zeta^{\l u^K, v(h_1)\r}\cdot \\
   \sum_{I\in \Sigma_K^{(n-k)}}\frac{w(I)}{|H_I|}
   &\sum_{h_2\in H_I}
   \prod_{i\in I}
    \phi(-\l u_i^I, zv-\tau v(h_1)+v(h_2)\r, \tau, \sigma),
\end{split} 
\end{equation}
where $u^K=\sum_{i\in K}u_i^K$.

\begin{note}
In the sum above with respect to $K\in\sigmk$ and $h_1\in \hat{H}_K$, 
the term corresponding to $K=*\in \Sigma^{(0)}$ and $h_1=0\in \hat{H}_*=0$ 
is equal to $\varv(\Delta,\V)$.
\end{note}

It is also sometimes useful to take a representative $v(h)$ of 
$h\in H_I$ such that
\begin{equation}\label{eqn:fIh}
 0\leq\l u_i^I,v(h)\r<1 \ \text{for all $i\in I$}. 
\end{equation}
Such a representative is unique. We denote the value 
$\l u_i^I,v(h)\r$ by $f_{I,h,i}$ for such a representative 
$v(h)$. If $h$ lies in $H_K$ for $K\in \sigmk$ contained in $I$, 
then $f_{I,h,i}=0$ for $i\not\in K$, and $f_{I,h,i}$ depends only 
on $K$ for $i\in K$ which we shall denote by $f_{K,h,i}$. 
The sum $\sum_{i\in K}f_{K,h,i}$ will be denoted by $f_{K,h}$. 
Note that \eqref{eqn:0} can be rewritten 
as
\begin{equation*}\label{eqn:0bis}
 \hat{H}_K=\{h\in H_K\mid f_{K,h,i}\not=0 
\quad \text{for any $i\in K$}. \}
\end{equation*} 

If we choose representatives $v(h)$ for $h\in H_I$ 
satisfying \eqref{eqn:fIh}, then 
\eqref{eq:hatvarv2} can be put in the form 
\begin{equation}\label{eq:hatvarv2bis}
 \begin{split}
 \hatvarv(\Delta,\V) 
 = \sum_{k=0}^n&\sum_{K\in\sigmk,h_1\in \hat{H}_K}
   \zeta^{f_{K,h}}\cdot \\
   \sum_{I\in \Sigma_K^{(n-k)}}\frac{w(I)}{|H_I|}
   &\sum_{h_2\in H_I}
   \prod_{i\in I}
    \phi(-\l u_i^I, zv-\tau v(h_1)+v(h_2)\r, \tau, \sigma). 
\end{split} 
\end{equation}

Let $N>1$ be an integer. For a rational number $f=\frac{s}{r}$ with 
$r$ relatively prime to $N$, we take an integer $d$ such that 
$dr\equiv 1 \bmod N$ and define 
\[ \breve{f}=ds. \]
The integer $\breve{f}$ is defined modulo $N$. 

Assume that $N$ is relatively prime to $|H_I|$ 
for all $I\in \sigmn$. We put $\sigma=\frac{k}{N}$ 
with $0<k<N$, and define the 
\emph{modified orbifold elliptic genus} 
$\brvarv(\Delta,\V)=\brvarv(\Delta,\V;z,\tau,\sigma)$ of 
\emph{level $N$} by 
\begin{equation*}\label{eq:brvarv}
 \begin{split}
 \brvarv(\Delta,\V) 
 = &\sum_{I\in\sigmn}
   \frac{w(I)}{|H_I|}\cdot \\
   &\sum_{(h_1,h_2)\in H_I\times H_I}
  \prod_{i\in I}\zeta^{\l u_i^I, v(h_1)\r\breve{}}
 \phi(-\l u_i^I, zv-\tau v(h_1)+v(h_2)\r, \tau, \sigma).
\end{split} 
\end{equation*}
It has also the following expression:
\begin{equation}\label{eq:brvarv2}
 \begin{split}
 \brvarv(\Delta,\V) 
 &= \sum_{k=0}^n\sum_{K\in\sigmk,h_1\in \hat{H}_K}
   \zeta^{\l u^K, v(h_1)\r\breve{}}\cdot \\
   &\sum_{I\in \Sigma_K^{(n-k)}}\frac{w(I)}{|H_I|}
   \sum_{h_2\in H_I}
    \prod_{i\in I}
    \phi(-\l u_i^I, zv-\tau v(h_1)+v(h_2)\r, \tau, \sigma) 
\end{split}
\end{equation}
If we choose representatives $v(h)$ satisfying \eqref{eqn:fIh}, 
then 
\begin{equation*}\label{eq:brvarv2bis}
\begin{split}
 \brvarv(\Delta,\V) 
 &= \sum_{k=0}^n\sum_{K\in\sigmk,h_1\in \hat{H}_K}
   \zeta^{\breve{f}_{K,h}}\cdot \\
   &\sum_{I\in \Sigma_K^{(n-k)}}\frac{w(I)}{|H_I|}
   \sum_{h_2\in H_I}
    \prod_{i\in I}
    \phi(-\l u_i^I, zv-\tau v(h_1)+v(h_2)\r, \tau, \sigma).
\end{split} 
\end{equation*}

\begin{prop}\label{prop:laurent}
Let $\varv(\Delta,\V) =
\sum_{s=0}^\infty \var_{s}(z) q^s$  
be the expansion
into power series, 
then $\zeta^{\frac{n}{2}}\var_{s}(z)$ belongs to 
$R(S^1)\otimes\Z[\zeta,\zeta^{-1}]$,
where $R(S^1)$ is identified with
$\Z[t,t^{-1}]$. Let $r$ be the least common multiple 
of $\{\vert H_I\vert\}_{I\in \sigmn}$. 
Then $\hatvarv(\Delta,\V)$ can be expanded in the form 
$\hatvarv(\Delta,\V) =
\sum_{s=0}^\infty \hatvar_{s}(z) q^s$, 
where $\zeta^{\frac{n}{2}}\hatvar_{s}(z)$ belongs to
$R(S^1)\otimes
\Z[\zeta^{\frac1{r}},\zeta^{-\frac1{r}}]$. Similarly 
$\brvarv(\Delta,\V)$ can be expanded in the form 
$\brvarv(\Delta,\V) =
\sum_{s=0}^\infty \brvar_{s}(z) q^s$, 
where $\zeta^{\frac{n}{2}}\brvar_{s}(z)$ belongs to
$R(S^1)\otimes
\Z[\zeta,\zeta^{-1}]/(\zeta^N)$.
\end{prop}
For the details of proof we refer to \cite{HM2}. 
We introduce an auxiliary variable $\tau_1$ with $\Im(\tau_1)>0$ 
and put $q_1=e^{2\pi\img\tau_1}$. 
The proof amounts to showing that 
\begin{equation*}\label{eq:laurent2}
\sum_{I\in \Sigma_K^{(n-k)}}\frac{w(I)}{|H_I|}
   \sum_{h_2\in H_I}
    \prod_{i\in I}
    \phi(-\l u_i^I, zv-\tau_1 v(h_1)+v(h_2)\r, \tau, \sigma)
\end{equation*}
is expanded in a power series
\[ \hatvar_{K,h_1}(z,\tau_1,\tau,\sigma)= 
 \sum_{s_1\in \Z,s_2\in \Z_{\geq 0}}\hatvar_{s_1,s_2}(z,q_1)\zeta^{s_1}q^{s_2}. \]
where $\zeta^{\frac{n}{2}}\hatvar_{s_1,s_2}(z,q_1)$ 
a finite sum of the expressions of the form 
\begin{equation}\label{eq:laurent3}
 \sum_{I\in\Sigma_K^{(n-k)}}\frac{w(I)}{|H_I|}
 \sum_{l\in \alpha_I}\sum_{h_2\in H_I}
 \frac{\chi_I(l,h_2)t^{\l l,v\r}q_1^{\l l,v(h_1)\r}}
 {\prod_{i\in I\setminus K}(1-\chi_I(u_i^I,h_2)^{-1}
 t^{\l -u_i^I,v\r}q_1^{\l -u_i^I,v(h_1)\r})}. 
\end{equation}
The expression \eqref{eq:laurent3} belongs to 
$\Z[t,t^{-1},q_1,q_1^{-1}]$ by Corollary \ref{coro:multiplicity}. 
From this fact the statement for $\hatvarv(\Delta,\V)$ follows. 
Note that $\hatvarv(\Delta,\V)$ does not have negative power of $q$ 
by \eqref{eq:hatvarv2bis}. 
The cases of $\varv(\Delta,\V)$ and $\brvarv(\Delta,\V)$ are similar. 

The elliptic genus $\var(\Delta,\V)=\var(\Delta,\V;\tau,\sigma)
\in (R(T)\otimes \Z[\zeta,\zeta^{-1}])[[q]]$, 
the orbifold elliptic genus $\hatvar(\Delta,\V)
=\hatvar(\Delta,\V;\tau,\sigma)\in 
(R(T)\otimes \Z[\zeta^{\frac1{r}},\zeta^{-\frac1{r}}])[[q]]$ 
and the modified orbifold elliptic genus 
$\brvar(\Delta,\V)
=\brvar(\Delta,\V;\tau,\sigma)\in 
(R(T)\otimes \Z[\zeta,\zeta^{-1}]/(\zeta^N))[[q]]$ 
of level $N$ are defined by
\[ v^*(\var(\Delta,\V))=\varv(\Delta,\V),\ 
  v^*(\hatvar(\Delta,\V))=\hatvarv(\Delta,\V)\ \ \text{and}
  \ \ v^*(\brvar(\Delta,\V))=\brvarv(\Delta,\V) \]
where one varies generic vectors $v$ in $L_{\V}$. 
\begin{note}
Once $\var(\Delta,\V)$, $\hatvar(\Delta,\V)$ and $\brvar(\Delta,\V)$ 
are defined as above we can define 
$\varv(\Delta,\V)$, $\hatvarv(\Delta,\V)$ and $\brvarv(\Delta,\V)$ 
by using the above formulas \emph{for any $v\in L$}. 
\end{note}

Let $N>1$ be an integer. When $\sigma=\frac{k}{N},\ 0<k<N$, the genera 
$\var(\Delta,\V)$ and $\hatvar(\Delta,\V)$ 
will be also called of \emph{level} $N$. 
\begin{theo}\label{theo:brrigid}
Let $(\Delta,\V)$ be a pair of complete simplicial multi-fan in a lattice 
$L$ of rank $n$ and 
a set of generating edge vectors. Let $N>1$ be an integer relatively 
prime to $|H_I|$ for every $I\in \sigmn$. If $c_1(\Delta,\V)$ is 
divisible by $N$, then the modified orbifold elliptic genus 
$\brvar(\Delta,\V;\tau,\sigma)$ of level $N$ constantly vanishes. 
\end{theo}

\begin{theo}\label{theo:hatrigid}
Let $(\Delta,\V)$ be a pair of complete simplicial multi-fan in a lattice 
$L$ of rank $n$ and 
a set of generating edge vectors. Let $N>1$ be an integer. If 
$c_1(\Delta,\V)$ is $T$-Cartier divisible by $N$, 
then the orbifold elliptic genus 
$\hatvar(\Delta,\V;\tau,\sigma)$ of level $N$ constantly vanishes. 
\end{theo}

\begin{theo}\label{theo:varrigid}
Let $(\Delta,\V)$ be a pair of complete simplicial multi-fan in a lattice 
$L$ of rank $n$ and 
a set of generating edge vectors. If 
$c_1(\Delta,\V)=0$, then the orbifold elliptic genus 
$\hatvar(\Delta,\V;\tau,\sigma)$ constantly vanishes. 
\end{theo}

The following examples are obtained by using Theorem 3.4 in 
\cite{HM2}. Consider the multi-fan $\Delta=(\Sigma,C,w^\pm)$ where 
\[ \sigmone=\{1,2,3\},\ \Sigma^{(2)}=\{\{1\},\{2\},\{3\}\}, \] 
and $C(1),C(2),C(3)$ are half lines in $\R^2$ generated by 
$e_1=(1,0),e_2=(0,1),-e_1-be_2$ respectively, where $b>0$ is an integer. 
We set $w^+(I)=1$ and $w^-(I)=0$ so that $w(I)=1$ for all 
$I\in \Sigma^{(2)}$. 

\begin{exam}
Define $\V=\{v_1,v_2,v_3\}$ by 
\[ v_1=e_1,v_2=e_2,v_3=-e_1-be_2. \] 
Then 
\[ \begin{split}\hatvar&(\Delta,\V)= \\
  &\left(\sum_{(m_1,m_2)\in \Z^2}t_1^{-m_1}t_2^{-m_2}
 \frac{(1-\zeta^{\frac{b+2}{b}})(1-\zeta^2q^{-bm_2})}
  {(1-\zeta q^{m_1})(1-\zeta q^{m_2})(1-\zeta q^{-m_1-bm_2})
 (1-\zeta^{\frac{2}{b}}q^{-m_2})}\right)\Phi(\sigma,\tau)^2.
 \end{split} \]

If $b$ is odd, then $c_1(\Delta,\V)$ is divisible by $b+2$, and 
we have
\[ \begin{split}\brvar&(\Delta,\V)= \\
  &\left(\sum_{(m_1,m_2)\in \Z^2}t_1^{-m_1}t_2^{-m_2}
  \frac{(1-\zeta^{b+2})(1-\zeta^2q^{-bm_2})}
  {(1-\zeta q^{m_1})(1-\zeta q^{m_2})(1-\zeta q^{-m_1-bm_2})
 (1-\zeta^{2l}q^{-m_2})}\right)\Phi(\sigma,\tau)^2=0, 
 \end{split} \]
where $b=2l-1$ and $\sigma=\dfrac{k}{b+2},\ 0<k<b+2$.
\end{exam}

\begin{exam}
Define $\V=\{v_1,v_2,v_3\}$ by 
\[ v_1=e_1,v_2=be_2,v_3=-e_1-be_2 \] 
so that $v_1+v_2+v_3=0$. 
Then 
\[ \begin{split}\hatvar&(\Delta,\V)= \\
  &\left(\sum_{(m_1,m_2)\in \Z^2}t_1^{-m_1}t_2^{-m_2}
  \frac{(1-\zeta^{\frac{3}{b}})(1-\zeta^2q^{-bm_2})}
  {(1-\zeta q^{m_1})(1-\zeta^{\frac{1}{b}}q^{m_2})(1-\zeta q^{-m_1-bm_2})
 (1-\zeta^{\frac{2}{b}}q^{-m_2})}\right)\Phi(\sigma,\tau)^2.
 \end{split} \]

If $b$ is relatively prime to $3$, then $c_1(\Delta,\V)$ is divisible 
by $3$, and we have
\[ \begin{split}&\brvar(\Delta,\V)= \\
  &\left(\sum_{(m_1,m_2)\in \Z^2}t_1^{-m_1}t_2^{-m_2}
  \frac{(1-\zeta^{3d})(1-\zeta^2q^{-bm_2})}
  {(1-\zeta q^{m_1})(1-\zeta^d q^{m_2})(1-\zeta q^{-m_1-bm_2})
 (1-\zeta^{2d}q^{-m_2})}\right)\Phi(\sigma,\tau)^2=0, 
 \end{split} \]
where $d$ is such that $db\equiv 1 \bmod \, 3$ and 
$\sigma=\dfrac{k}{3},\ k=1,2$.
\end{exam}

\begin{exam}
If $b=2$ in Example 1, then $c_1(\Delta,\V)$ is $T$-Cartier divisible 
by $2$. This corresponds to the toric variety $\Proj^2(2,1,1)$, see 
Corollary \ref{coro:N=n}. In this case 
\[ \hatvar(\Delta,\V)=
  \left(\sum_{(m_1,m_2)\in \Z^2}t_1^{-m_1}t_2^{-m_2}
  \frac{(1-\zeta^2)(1+\zeta q^{-m_2})}
  {(1-\zeta q^{m_1})(1-\zeta q^{m_2})(1-\zeta q^{-m_1-2m_2})}
 \right)\Phi(\sigma,\tau)^2.
  \]
\end{exam}

\begin{rema}
For ordinary fans $c_1(\Delta,\V)$ has infinite order when the vectors
$v_i$ are taken primitive. But 
there are examples of general multi-fans with vanishing 
$c_1(\Delta,\V)$. 
\end{rema}

Proofs of Theorem \ref{theo:brrigid}, Theorem \ref{theo:hatrigid} 
and Theorem \ref{theo:varrigid} will be given in Section 4. 
\par 
\vspace*{0.3cm}
The elliptic genus $\var(\Delta,\V)$ reduces to the 
so-called $T_y$-genus for $q=0$ if it is multiplied by $\zeta^{n/2}$ 
and if $\zeta$ is substituted by $-y$. 
Namely 
\[ T_y(\Delta,\V)=\sum_{I\in \sigmn}\frac{w(I)}{|H_I|}
   \sum_{h\in H_I}\prod_{i\in I}
   \frac{1+y\chi_I(u_i^I,h)^{-1}t^{-u_i^I}}
 {1-\chi_I(u_i^I,h)^{-1}t^{-u_i^I}}. \]

In \cite{HM2} it was shown that $T_y$-genus $T_y(\Delta,\V)$ 
had the following expression. 
\begin{equation}\label{eq:h}
 T_y(\Delta,\V)=\sum_{k=0}^nh_k(\Delta)(-y)^k, 
\end{equation}
where $h_k(\Delta)$ is defined by
\[ h_k(\Delta)=\sum_{I\in \sigmn,\ \mu(I)=k}w(I)\quad 
\text{with}\ \
 \mu(I)=\#\{i\in I\mid \l u_i^I,v\r >0\}. \]
Here $v$ is a generic vector but $h_k(\Delta)$ does not
depend on the choice of $v$. 
The equality \eqref{eq:h} shows that $T_y(\Delta,\V)$ is rigid 
and is independent of $\V$. So we write it simply $T_y(\Delta)$. 

In \cite{HM2} it was also shown that $T_y(\Delta)$ 
had the following expression. 
\begin{equation}\label{eq:T_y}
 T_y(\Delta)=\sum_{k=0}^ne_k(\Delta)(-1-y)^{n-k}
\end{equation}
where $e_k(\Delta)=\sum_{J\in \sigmk}\deg(\Delta_J)$.

\begin{note}
We have $h_k(\Delta)=h_{n-k}(\Delta)$. $h_0(\Delta)=T_0(\Delta)$ 
is the Todd genus of $\Delta$, 
and $h_n(\Delta)=\deg(\Delta)$ by definition of the latter. 
Hence $T_0[\Delta]$ equals $\deg(\Delta)$, cf. \cite{HM2}.
\end{note}

We define \emph{orbifold $T_y$-genus} of $(\Delta,\V)$ 
by
\begin{equation*}\label{eq:hatTy}
 \hatT_y(\Delta,\V)=\sum_{k=0}^n\sum_{K\in\sigmk}\sum_{h\in\hatH_K}
 (-y)^{f_{K,h}}T_y(\Delta_K).
\end{equation*}
It is equal to the degree zero term of 
$\zeta^{n/2}\hatvar(\Delta,\V)$ with $\zeta$ substituted by $-y$. 

When $|H_I|$ is relatively prime to an integer $N>1$ for every 
$I\in \sigmn$, the \emph{modified orbifold $T_y$-genus} of level $N$ 
of $(\Delta,\V)$ is defined by 
\begin{equation*}\label{eq:brTy}
\brT_y(\Delta,V)=\sum_{k=0}^n\sum_{K\in\sigmk}\sum_{h\in\hatH_K}
 (-y)^{\breve{f}_{K,h}}T_y(\Delta_K),
\end{equation*}
where $-y=e^{2\pi\img\frac{l}{N}},\ 0<l<N$. 
It is equal to the degree zero term of 
$\zeta^{n/2}\brvar(\Delta,\V)$ with $\zeta$ substituted by $-y$ . 
 
Suppose that $c_1(\Delta,\V)$ is $T$-Cartier. 
Then $\l u^I,v(h)\r=\l \iota_I^*(c_1(\Delta,\V)), v(h)\r$ is 
an integer 
because $\iota_I^*(c_1(\Delta,\V))$ lies in $L^*$. Since 
$f_{K,h}\equiv \l\iota_I^*(c_1(\Delta,\V)), v(h)\r \bmod \Z$ for 
$I\supset K$, $f_{K,h}$ is an integer for any $K\in \sigmk$ and 
$h\in \hatH_K$. It follows that $\hatT_y(\Delta,\V)$ is a 
polynomial in $-y$. This is the case in particular when 
$c_1(\Delta,\V)$ is $T$-Cartier divisible by $N$. 

\begin{prop}\label{prop:hatTdiv}
Under the situation of Theorem \ref{theo:hatrigid}, the 
orbifold $T_y$-genus $\hatT_y(\Delta,\V)$ is a polynomial with integer 
coefficients divisible by \[ \sum_{k=0}^{N-1}(-y)^k .\]
\end{prop}
In fact, by Theorem \ref{theo:hatrigid}, $\hatvar(\Delta,\V)$ vanishes 
for $-y=e^{2\pi\img \frac{k}{N}},\ 0<k<N,$ under the assumption, and 
hence its degree $0$ term $\hatT_y(\Delta,\V)$ does so. 
Thus it is a polynomial in $-y$ which vanishes 
for $-y=e^{2\pi\img \frac{k}{N}},\ 0<k<N$. Hence $\hatT_y(\Delta,\V)$ 
must be divisible by $\sum_{k=0}^{N-1}(-y)^k$. 

The following two propositions are corollaries of 
Theorem \ref{theo:brrigid} 
and Theorem \ref{theo:varrigid}. 
\begin{prop}\label{prop:brTdiv}
Under the situation of Theorem \ref{theo:brrigid}, the modified 
orbifold $T_y$-genus $\brT_y(\Delta,\V)$ of level $N$ vanishes. 
\end{prop}

\begin{prop}\label{prop:hatTvanish}
Under the situation of Theorem \ref{theo:varrigid}, the 
orbifold $T_y$-genus $\hatT_y(\Delta,\V)$ vanishes. 
\end{prop}

\section{Proofs of main theorems}\label{sec:proof}
Let $\Delta$ be a complete simplicial multi-fan and 
$\V$ a set of prescribed edge vectors. Let $v\in L_\V$ a generic vector. 
We put $H_I=L/L_{I,\V}$ for $I\in \sigmn$ as before. For $A\in SL_2(\Z)$ 
we set 
\[ (\hatvarv)^A(\Delta,\V;z,\tau,\sigma)=
 \hatvarv(\Delta,\V;A(z,\tau),\sigma). \]
Similarly if $N>1$ is an integer relatively prime to $|H_I|$ for 
all $I\in \sigmn$, then we set 
\[ (\brvarv)^A(\Delta,\V;z,\tau,\sigma)=
 \brvarv(\Delta,\V;A(z,\tau),\sigma),\ \sigma=\frac{k}{N}, \ 0<k<N. \]

If $c_1(\Delta,\V)$ is divisible by $N$, then  
it follows from \eqref{eq:Ndivide} that
\begin{equation}\label{eq:equilib}
 \l u^I,v\r=\l \iota_I^*(c_1^T(\Delta,\V)),v\r =
 N\l\iota_I^*(x),v\r+\l u,v\r .
\end{equation}
for any $I\in \sigmn$. In particular the $\bmod\, N$ value of 
the integer $\l u^I,v\r$ is equal to the $\bmod\, N$ value of $\l u,v\r$ 
and is independent of $I$. It will be denoted by $h(v)$. When 
$x$ is $T$-Cartier in \eqref{eq:Ndivide} $\l u^I,v\r$ and 
$\l\iota_I^*(x),v\r$ are integers for any vector $v\in L$. 

\begin{lemm}\label{lemm:hv}
Assume that $c_1^T(\Delta,\V)$ is divisible by $N$. 
If one of the following conditions is satisfied, then there exists 
a generic vector $v\in L$ such that $h(v)$ is defined and relatively 
prime to $N$.
\begin{enumerate}
\renewcommand{\labelenumi}{(\alph{enumi})} 
\item $N$ is relatively prime to $|H_I|$ for all $I\in \sigmn$.
\item $c_1(\Delta,\V)$ is $T$-Cartier divisible by $N$. 
\item $L_\V=L_{I,\V}$ for all $I\in \sigmn$. 
\end{enumerate}
Such a vector $v$ is taken in $L_\V$ for the cases (a) and (c). 
\end{lemm}
\begin{proof}
Fix an element $I\in \sigmn$. Since the $\u$ form a basis of 
$L_{I,\V}^*$, there is a $v\in L_{I,\V}$ 
such that $\l u^I,v\r$ takes a given integer value. In particular there is a 
$v\in L_{I,\V}$ such that $\l u^I,v\r$ is relatively prime to $N$. 

In the case (a) the index of $L_\V$ in $L_{I,\V}$ is relatively
prime to $N$ as is easily seen. Hence $v$ as above can be taken 
in $L_\V$. In the case (c) $v$ lies in $L_\V=L_{I,\V}$. 

Since $c_1^T(\Delta,\V)$ is divisible by $N$, the value 
$\l u^I,v\r \ \bmod\, N$ is independent of $I$ and equal to $h(v)$. 
\end{proof}

\begin{lemm}\label{lemm:brA}
Assume that $N$ is relatively prime to $|H_I|$ for all $I\in \sigmn$ 
and $c_1(\Delta,\V)$ is divisible by $N$. Then 
$(\brvarv)^A(\Delta,\V;z,\tau,\sigma)$ with $\sigma=\frac{k}{N},\ 0<k<N$ 
has the following expression. 
\begin{equation}\label{eq:brA}
 \begin{split}
 (\brvarv)^A&(\Delta,\V;z,\tau,\sigma)=
 e^{\pi\img(nc(c\tau+d)\sigma^2-2c\l u,v\r z\sigma)} \\
 &\sum_{I\in \sigmn}\frac{w(I)}{|H_I|}\sum_{(h_1,h_2)\in H_I\times H_I}
 e^{-2\pi\img\l \iota_I^*(x),dkv(h_1)\r}
 e^{-2\pi\img\l \iota_I^*(x), ck(zv+v(h_2))\r} \\
 &\quad \prod_{i\in I}e^{2\pi\img\l u_i^I,(c\tau+d)\sigma v(h_1)\r}
 \phi(-\l u_i^I, zv-\tau v(h_1)+v(h_2)\r,\tau,(c\tau+d)\sigma).
 \end{split}
\end{equation}
\end{lemm}

\begin{lemm}\label{lemm:hatA}
Assume that 
$c_1(\Delta,\V)$ is $T$-Cartier divisible by $N$. 
Then $(\hatvarv)^A(\Delta,\V;z,\tau,\sigma)$ with $\sigma=\frac{k}{N},\ 0<k<N,$ 
has the following expression. 
\begin{equation}\label{eq:hatA}
 \begin{split}
 (\hatvarv)^A(&\Delta,\V;z,\tau,\sigma)=
 e^{\pi\img(nc(c\tau+d)\sigma^2-2c\l u,v\r z\sigma)} \\
 &\sum_{I\in \sigmn}\frac{w(I)}{|H_I|}\sum_{(h_1,h_2)\in H_I\times H_I}
 e^{-2\pi\img\l \iota_I^*(x), ckzv\r} \\
 &\quad \prod_{i\in I}e^{2\pi\img\l u_i^I,(c\tau+d)\sigma v(h_1)\r}
 \phi(-\l u_i^I, zv-\tau v(h_1)+v(h_2)\r,\tau,(c\tau+d)\sigma).
 \end{split}
\end{equation}
\end{lemm}

\begin{proof}
We first prove Lemma \ref{lemm:brA}. By definition we have 
\[ 
\begin{split}
 (\brvarv)&^A(\Delta,\V;z,\tau,\sigma) 
 = \sum_{I\in\sigmn}
   \frac{w(I)}{|H_I|}\cdot \\
   &\sum_{(h_1,h_2)\in H_I\times H_I}
  \prod_{i\in I}\zeta^{\l u_i^I, v(h_1)\r\breve{}}
 \phi(-\l u_i^I, \frac{zv-(a\tau+b)v(h_1)+(c\tau+d)v(h_2)}{c\tau+d}\r,
 A\tau, \sigma).
\end{split} \] 
Using \eqref{eq:transformula1} 
we get
\begin{equation}\label{eq:brA2}
 \begin{split}
 \prod_{i\in I}\zeta^{\l u_i^I, v(h_1)\r\breve{}}
 &\phi(-\l u_i^I, \frac{zv-(a\tau+b)v(h_1)+(c\tau+d)v(h_2)}{c\tau+d}\r,
 A\tau, \sigma)= \\
 &\zeta^{\sum_i\l\u,v(h_1)\r\breve{}}
 e^{\pi\img\left(nc(c\tau+d)\sigma^2+2c\l u^I,-zv+(a\tau+b)v(h_1)
 -(c\tau+d)v(h_2)\r\sigma\right)} \\
 & \quad \prod_{i\in I}\phi(-\l\u,zv-(a\tau+b)v(h_1)+(c\tau+d)v(h_2)\r,
 \tau,(c\tau+d)\sigma).
 \end{split} 
\end{equation}

We have
\[ c\left((a\tau+b)v(h_1)-(c\tau+d)v(h_2)\right)=
 -v(h_1)+(c\tau+d)(av(h_1)-cv(h_2)). \] 
Noting that $\breve{f}_1+\breve{f_2}\equiv (f_1+f_2)\breve{}\ \bmod N$ 
we see that 
\[ \zeta^{\sum_i\l\u,v(h_1)\r\breve{}}=\zeta^{\l u^I,v(h_1)\r\breve{}}.
 \]
Since $c_1(\Delta,\V)$ is divisible by $N$, we get from \eqref{eq:Ndivide} 
\[ \l u^I,v(h_1)\r =N\l\iota_I^*(x),v(h_1)\r+\l u,v(h_1)\r, \]
and hence $ \l u^I,v(h_1)\r\breve{}\equiv \l u,v(h_1)\r \ \bmod N$ 
because $u\in L^*$ and $\l u,v(h_1)\r\in \Z$. 
Therefore 
\[ \zeta^{\sum_i\l\u,v(h_1)\r\breve{}}=\zeta^{\l u,v(h_1)\r}
=e^{2\pi\img\l u,v(h_1)\r\sigma}, \]
and 
\begin{equation}\label{eq:4-1}
 \zeta^{\sum_i\l\u,v(h_1)\r\breve{}}e^{-2\pi\img \l u^I,v(h_1)\r\sigma}
 =e^{-2\pi\img k\l\iota_I^*(x),v(h_1)\r}. 
\end{equation}
 
Similarly we have by \eqref{eq:equilib}
\begin{equation}\label{eq:4-2}
 e^{2\pi\img\l u^I, zv\r\sigma}=
 e^{2\pi\img\l u,zv\r\sigma}e^{2\pi\img k\l\iota_I^*(x),zv\r}. 
\end{equation} 

Let $\rho:H_I\times H_I\to H_I\times H_I$ be the map defined by 
\[ \rho(h_1,h_2)=(\barh_1,\barh_2)=(ah_1-ch_2, -bh_1+dh_2). \]
$\rho$ is bijective and its inverse is given by 
\[ \rho^{-1}(\barh_1,\barh_2)=(d\barh_1+c\barh_2,b\barh_1+a\barh_2). \]
Then $av(h_1)-cv(h_2)$ and $-bv(h_1)+dv(h_2)$ are representatives of 
$\barh_1$ and $\barh_2$ which we shall denote by $v(\barh_1)$ and 
$v(\barh_2)$ respectively. We then have 
\begin{equation}\label{eq:4-3}
 v(h_1)=dv(\barh_1)+cv(\barh_2). 
\end{equation}

In view of \eqref{eq:4-1}, \eqref{eq:4-2} and \eqref{eq:4-3} 
the right hand side of \eqref{eq:brA2} is equal to
\begin{equation}\label{eq:brA3}
 \begin{split}
 &e^{\pi\img\left(nc(c\tau+d)\sigma^2-2c\l u,v\r z\sigma\right)}
 e^{-2\pi\img\l\iota_I^*(x),dkv(\barh_1)\r}
 e^{-2\pi\img\l\iota_I^*(x),ck(zv+v(\barh_2))\r}
 \\
 & \quad e^{2\pi\img\l u^I,(c\tau+d)\sigma v(\barh_1)\r}
 \prod_{i\in I}\phi(-\l\u,zv-v(\barh_1)\tau+v(\barh_2)\r,
 \tau,(c\tau+d)\sigma).
 \end{split} 
\end{equation}
Summing up over $(h_1,h_2)$ is the same as summing up over 
$(\barh_1,\barh_2)$. Hence from \eqref{eq:brA3} we get 
\eqref{eq:brA} with $h_i$ replaced by $\barh_i$ for $i=1,2$. 
This proves Lemma \ref{lemm:brA}. 

As to Lemma \ref{lemm:hatA} we have 
\[ \zeta^{\sum_i\l\u,v(h_1)\r}
  e^{-2\pi\img \l u^I,v(h_1)\r\sigma}=1 \]
instead of \eqref{eq:4-1}. The rest of the proof is entirely similar 
to that of Lemma \ref{lemm:brA}. 
\end{proof}

\begin{lemm}\label{lemm:brnopole}
Assume that $N$ is relatively prime to $|H_I|$ for all $I\in \sigmn$ 
and $c_1(\Delta,\V)$ is divisible by $N$. Then the meromorphic 
function $(\brvarv)^A(\Delta,\V;z,\tau,\sigma)$ in 
$z$ and $\tau$ with $\sigma=\frac{k}{N},\ 0<k<N$, has no pole at $z\in \R$. 
\end{lemm}

\begin{lemm}\label{lemm:hatnopole}
Assume that $c_1(\Delta,\V)$ is $T$-Cartier divisible by $N$. 
Then the meromorphic 
function $(\hatvarv)^A(\Delta,\V;z,\tau,\sigma)$ in 
$z$ and $\tau$ with $\sigma=\frac{k}{N},\ 0<k<N$, has no pole at $z\in \R$. 
\end{lemm}

\begin{proof}
The expression \eqref{eq:brA} in Lemma \ref{lemm:brA} of the 
function $(\brvarv)^A(\Delta,\V;z,\tau,\sigma)$ can be rewritten 
in the following form as can be seen in a similar way to 
\eqref{eq:brvarv2}. 
\begin{equation*}
 \begin{split}
 &(\brvarv)^A(\Delta,\V;z,\tau,\sigma) \\ 
  &=e^{\pi\img(nc(c\tau+d)\sigma^2-2c\l u,v\r z\sigma)}
   \sum_{k=0}^n\sum_{K\in\sigmk,h_1\in \hat{H}_K}
   e^{-2\pi\img\l u^K, (c\tau+d)\sigma v(h_1)\r}\cdot \\
   &\sum_{I\in \Sigma_K^{(n-k)}}\frac{w(I)}{|H_I|}
   e^{-2\pi\img\l \iota_I^*(x), ck(zv+v(h_2))\r}
   \sum_{h_2\in H_I}
    \prod_{i\in I}
    \phi(-\l u_i^I, zv-\tau v(h_1)+v(h_2)\r, \tau, (c\tau+d)\sigma).
\end{split} 
\end{equation*}
Hence, in order to prove Lemma \ref{lemm:brnopole}, it is 
sufficient to prove that 
\begin{equation*}
\sum_{I\in \Sigma_K^{(n-k)}}\frac{w(I)}{|H_I|}
   e^{-2\pi\img\l \iota_I^*(x), ck(zv+v(h_2))\r}
   \sum_{h_2\in H_I}
    \prod_{i\in I}
    \phi(-\l u_i^I, zv-\tau v(h_1)+v(h_2)\r, \tau, (c\tau+d)\sigma), 
\end{equation*}
or, replacing $ck$ by $m$ and $(c\tau+d)\sigma$ by $\sigma$, 
\begin{equation}\label{eq:brvarv3}
\sum_{I\in \Sigma_K^{(n-k)}}\frac{w(I)}{|H_I|}
   e^{-2\pi\img\l \iota_I^*(x), m(zv+v(h_2))\r}
   \sum_{h_2\in H_I}
    \prod_{i\in I}
    \phi(-\l u_i^I, zv-\tau v(h_1)+v(h_2)\r, \tau,\sigma), 
\end{equation}
has no pole at $z\in \R$
for any fixed $K\in \sigmk$ and $h_1\in \hatH_K$. 
Note that 
\[ e^{-2\pi\img\l \iota_I^*(x), m(zv+v(h_2))\r}=
 e^{-2\pi\img\l u^K,m\tau v(h_1)\r}
 e^{-2\pi\img\l \iota_I^*(x),m(zv-\tau v(h_1)+v(h_2))\r}. \]
By a similar argument to the proof of Proposition \ref{prop:laurent},
we see that \eqref{eq:brvarv3} can be expanded 
in the form 
\[ e^{-2\pi\img\l u^K,m\tau v(h_1)\r}
 \sum_{s=0}^\infty (\brvar)^A_{K,h_1,s}(z) q^s, \]
where $\zeta^{\frac{n}{2}}(\brvar)^A_{K,h_1,s}(z)$ belongs to
$R(S^1)\otimes
\Z[\zeta^{\frac1{r}},\zeta^{-\frac1{r}}]/(\zeta^N)$. 
From this we can conclude that \eqref{eq:brvarv3} has no pole at $z\in \R$. 
We refer to Lemma in Section 7 of \cite{Hir}. See also Section 5 of 
\cite{HM2}. This finishes the proof of Lemma \ref{lemm:brnopole}. 

The proof of Lemma \ref{lemm:hatnopole} is similar. 
We use Lemma \ref{lemm:hatA} instead of Lemma \ref{lemm:brA}. 
Then we have only to note 
that 
\[ e^{2\pi\img \l \iota_I^*(x),ckzv\r}=
  e^{2\pi\img \l \iota_I^*(x),ck(zv+v(h_2))\r} \]
which holds because $\iota_I^*(x)\in L$ and hence 
$\l\iota_I^*(x), v(h_i)\r\in \Z,\ i=1,2$. 
\end{proof}

We now proceed to the proof of Theorem \ref{theo:brrigid}. We follow 
\cite{L} for the idea of proof.  We first 
show that $\brvarv(\Delta,\V;z,\tau,\sigma)$ is a constant as 
a function of $z$. 

We regard 
$\brvarv(\Delta,\V;z,\tau,\sigma)$ as a meromorphic function of $z$. 
By the transformation law \eqref{eq:transformula1} $\phi(z,\tau,\sigma)$ 
is an elliptic function in $z$ with respect to the lattice 
$\Z\cdot N\tau\oplus \Z$ for $\sigma=\frac{k}{N}$ with $0<k<N$. 
Hence $\brvarv(\Delta,\V;z,\tau,\sigma)$ of level $N$ is also 
an elliptic function 
in $z$. Thus, in order to show that $\brvarv(\Delta,\V;z,\tau,\sigma)$ is a 
constant function it suffices to show that it does not have poles. 

Assume that $z$ 
is a pole. Then $1-t^mq^r\alpha =0$ for some integer $m\not=0$, some
rational number $r$ and a root of unity $\alpha$. 
Consequently there are intergers $m_1\not=0$ and $k_1$ such that
$ m_1z+k_1\tau \in\R$. Then there is an element 
$A= \begin{pmatrix} a & b \\ c & d \end{pmatrix} \in SL_2(\Z)$ such 
that
\[ \frac{z}{c\tau+d}\in \R. \]
Since
\[ \brvarv(\Delta,\V;z,\tau,\sigma)=
 \brvarv(\Delta,\V;A^{-1}(\frac{z}{c\tau+d},A\tau),\sigma)=
 (\brvarv)^{A^{-1}}(\Delta,\V;\frac{z}{c\tau+d},A\tau,\sigma), \]
$(\brvarv)^{A^{-1}}(\Delta,\V;\frac{z}{c\tau+d},A\tau,\sigma)$ 
must have a pole 
at $\frac{z}{c\tau+d}\in \R$. But this contradicts Lemma 
\ref{lemm:brnopole}. This contradiction proves that 
$\brvarv(\Delta,\V;z,\tau,\sigma)$ 
can not have a pole. 

Since $\brvarv(\Delta,\V;z,\tau,\sigma)$ is a constnat function in $z$ 
for every generic vector $v\in L$, the equivariant modified orbifold 
genus $\brvar(\Delta,\V;\tau,\sigma)$ is constant as a function 
on $T$. That constant is equal to $\brvarv(\Delta,\V;z,\tau,\sigma)$ for 
any $v$. 
On the other hand, using \eqref{eq:transformula1} and the fact that 
$h(v) \bmod N$ is independent of $I\in \sigmn$, we have 
\begin{equation}\label{eq:hv}
 \brvarv(\Delta,\V;\tau,\sigma)=\brvarv(\Delta,\V;z+\tau,\tau,\sigma)=
 \zeta^{h(v)}\brvarv(\Delta,\V;z,\tau,\sigma)=\brvarv(\Delta,\V;\tau,\sigma). 
\end{equation}
We choose a generic vector $v$ such that $h(v)$ is relatively 
prime to $N$. Such a $v$ exists by Lemma \ref{lemm:hv}. Then 
$\zeta^{h(v)}$ is not equal to $1$. Hence 
\eqref{eq:hv} implies that the constant $\brvar(\Delta,\V;\tau,\sigma)$ 
must vanish. This finishes the proof of Theorem \ref{theo:brrigid}. 

The proof of Theorem \ref{theo:hatrigid} is similar. We have only 
to use Lemma 
\ref{lemm:hatnopole} instead of Lemma \ref{lemm:brnopole}. 

For the proof of Theorem \ref{theo:varrigid} we first note 
that the condition $c_1(\Delta,\V)=0$ 
is equivalent to $u=c_1^T(\Delta,\V)\in L^*$. In particular 
$\l u^I,v\r =\l u,v\r\in \Z$ for all $I\in\sigmn$. Then we have the equality 
\begin{equation*}
 \begin{split}
 (\hatvarv)^A&(\Delta,\V;z,\tau,\sigma)=
 e^{\pi\img(nc(c\tau+d)\sigma^2-2c\l u,v\r z\sigma)}
 \sum_{I\in \sigmn}\frac{w(I)}{|H_I|}\cdot \\
 &\sum_{(h_1,h_2)\in H_I\times H_I}
 \prod_{i\in I}e^{2\pi\img\l u_i^I,(c\tau+d)\sigma v(h_1)\r}
 \phi(-\l u_i^I, zv-\tau v(h_1)+v(h_2)\r,\tau,(c\tau+d)\sigma)
 \end{split}
\end{equation*}
which can be proved in a similar way to \eqref{eq:hatA}. 
From this equality we can conclude that 
$(\hatvarv)^A(\Delta,\V;z,\tau,\sigma)$ has no pole at $z\in \R$ 
as in the proof of Lemma \ref{lemm:hatnopole}. 
Let $N$ be an integer greater than $1$. Then one sees, as in the proof of 
Theorem \ref{theo:brrigid}, that 
$\hatvarv(\Delta,\V;z,\tau,\sigma)$ is a constant for 
$\sigma=\frac{k}{N},\ 0<k<N$. Since this is true for any integer $N>1$ 
and for any generic $v\in L_\V$, 
$\hatvarv(\Delta,\V;z,\tau,\sigma)$ must be a constant equal to 
$\hatvar(\Delta,\V;\tau,\sigma)$. Moreover we have 
\[  \hatvarv(\Delta,\V;z+\tau,\tau,\sigma)=
 \zeta^{\l u,v\r}\hatvarv(\Delta,\V;z,\tau,\sigma). \]
Then choosing 
$v\in L_\V$ such that $h(v)=\l u,v\r\not=0$ we see that 
$\hatvar(\Delta,\V;\tau,\sigma)$ must be equal to $0$.

\section{Applications}
Let $\Delta=(\Sigma,C,\pm)$ be a complete simplicial multi-fan 
in a lattice $L$ of rank $n$ . 
In this section a vector $v_i\in L$ generating the cone $C(i)$ 
for each $i\in \sigmone$ will always be taken primitive 
so that $\V=\{v_i\}$ 
is determined by $\Delta$. Thus $T_y(\Delta,\V)$ and 
$\hatT_y(\Delta,\V)$ will be simply written $T_y(\Delta)$ and 
$\hatT_y(\Delta)$. Similarly $c_1(\Delta,\V)$ is written 
$c_1(\Delta)$.

Also the following condition will be assumed throughout this section. 
\begin{equation}\label{eq:Toddone}
\text{\emph{The Todd genus $T_0(\Delta)$ is equal to $1$ and $w(I)=1$ 
for all $I\in \sigmn$.}}
\end{equation}
\begin{note}
The condition \eqref{eq:Toddone} is always satisfied by 
complete simplicial ordinary fans. See e.g. \cite{Ful}.
\end{note}

\begin{lemm}\label{lemm:Toddone1}
Under the condition \eqref{eq:Toddone} we have $T_0(\Delta_K)=1$ for all 
$K\in \sigmk$ and $w_K(I)=1$ for all $I\in \Sigma_K^{(n-k)}$ and 
$K\in \sigmk$. 
\end{lemm}

\begin{proof}
Take $I\in \Sigma_K^{n-k}$. $I$ is an element in $\sigmn$ such that 
$K\subset I$ and $w_K(I)=w(I)=1$ by definition and by assumption. 

As was remarked in Note after 
\eqref{eq:T_y}, $T_0(\Delta)$ is equal to $\deg(\Delta)$. We shall show that
$\deg(\Delta_K)=1$ for each $K\in \sigmk$ which will prove Lemma. 
Take a generic vector $\bar{v}$ in $L^K_\R$ and a 
generic vector $v$ in $L_\R$ which projects into $\bar{v}$. 
Since $\deg(\Delta)=1$ and $w(I)=1$ for all $I\in \sigmn$, there is 
a unique $I\in \sigmn$ with $K\subset I$ such that $v$ is contained 
in $C(I)$. Then $C_K(I)$ contains $\bar{v}$. 

We may assume that $v$ is chosen in such a way that $\l u_i^I,v\r$ 
is a sufficiently large positive number for each $i\in K$. 
Assume that there is another $I'\in \Sigma_K^{n-k}$ such that 
$C_K(I')$ contains $\bar{v}$. Then, $\l u_i^{I'}, v_i\r >0$ for all 
$i\in I'\setminus K$. From the fact that $\l u_i^I,v\r$ is 
sufficiently large for every $i\in K$, it follows that 
$\l u_i^{I'},v\r>0$ also for all
$i\in K$ and hence $v\in C(I')$. This contradicts the fact 
that $I$ is the unique element in $\sigmn$ such that 
$I\supset K$ and $v\in C(I)$. Hence $I\in \Sigma_K^{n-k}$ is the 
unique elemnet which contains $\bar{v}$. Since $w_K(I)=1$ we have
$\deg(\Delta_K)=1$. 
\end{proof}

\begin{lemm}\label{lemm:Toddone2}
Under the condition \eqref{eq:Toddone} the following equalities hold. 
\begin{equation}\label{eq:hq}
 h_k(\Delta)=\#\{I\in \sigmn \mid \mu(I)=k\},
\end{equation}
in \eqref{eq:h} and 
\begin{equation}\label{eq:ek}
  e_k(\Delta)=\#\sigmk 
\end{equation}
in (\ref{eq:T_y}).
\end{lemm}
\begin{proof}
\eqref{eq:hq} is immediate. 
$\deg(\Delta_K)=1$ for all $K\in \Sigma$ by Lemma \ref{lemm:Toddone1}  
Then \eqref{eq:ek} follows.
\end{proof}

\begin{lemm}\label{lemm:hpositive}
$h_k(\Delta)>0$ for $0\leq k\leq n$.
\end{lemm}
\begin{proof}
Let $v$ be a generice vector. For $I\in \sigmn$ we put 
\[ \mu_v(I)=\#\{i\in I\mid \l \u ,v\r >0\}. \]
Then $h_k(\Delta)=\#\{I\in \sigmn \mid \mu_v(I)=k\}$ 
by \eqref{eq:hq}. Fixing $I$ we put
\[ \sigma_J=\{v\in L_\R \mid \l u_j^I, v\r>0 \ \text{for}\ 
j\in J,\ \l u_i^I, v\r <0 \ \text{for}\ i\not\in J \} \]
for $J\subset I$. 
Then the collection $\{\sigma_J\}_{J\subset I}$ decomposes 
$L_\R\setminus \{v\mid \l\u,v\r=0\ \ \text{for some $i\in I$}\}$ 
into connected components. 
If we take $J\subset I$ in $\sigmk$ and $v$ in $\sigma_J$, then 
$\mu_v(I)=k$. This proves that $h_k(\Delta)>0$. 
\end{proof}
\begin{note}
It is known that $h_{k-1}(\Delta)\leq h_k(\Delta)$ for 
$0\leq k\leq[\frac{n}{2}]$ 
for the fan associated to a complete $\Q$-factorial projective 
toric variety, cf. \cite{Ful}.
\end{note}

We shall also use the following fact.
\begin{equation}\label{eq:fpositive}
 0<f_{K,h}<k \ \ \text{\emph{for}}\ K\in \sigmk,\ k>0, \ 
 \text{\emph{and}}\ h\in \hatH_K.  
\end{equation}
In fact $\l u_i^K, v(h)\r\not\in \Z$ for $h\in \hatH_K$ and 
$i\in K$ by \eqref{eqn:0}. Hence $0<f_{K,h,i}<1$ and 
$0<f_{K,h}=\sum_{i\in K}f_{K,h,i}<k$. 

\begin{prop}\label{prop:N=n+1}
Let $\Delta$ be a complete multi-fan of dimension $n$ 
satisfying condition \eqref{eq:Toddone}. If 
$c_1(\Delta)$ is $T$-Cartier divisible by an integer $N>1$, then 
$N$ is equal to or less than $n+1$. In the extremal case $N=n+1$ 
the multi-fan $\Delta$ is non-singular and
the $T_y$-genus must be of the form
\begin{equation}\label{eq:Tyn+1}
 T_y(\Delta)=\sum_{k=0}^n(-y)^k.
\end{equation} 
\end{prop}
\begin{proof}
Suppose that $c_1(\Delta)$ is $T$-Cartier divisible by $N$. 
Then, 
by Propsition \ref{prop:hatTdiv}, $\hatT_y(\Delta)$ is 
divisible by $\sum_{k=0}^{N-1}(-y)^k$.
On the other hand it is a polynomial of degree $n$ with constant
term $T_0(\Delta)=1$. Therefore we must have
$N-1\leq n$.

Suppose that $N=n+1$. Then
the same reasoning as above shows that 
\[ \hatT_y(\Delta)=\sum_{k=0}^{n}(-y)^k. \]
On the other hand we have 
\[ T_y(\Delta)=\sum_{k=0}^{n}h_k(-y)^k \]
with $h_k>0$ by Lemma \ref{lemm:hpositive}. 
Hence we must have $T_y(\Delta)=\hatT_y(\Delta)$. It also 
follows that $\Delta$ is non-singular. For otherwise 
there would be an extra term $(-y)^{f_{K,h}}T_y(\Delta_K)$ 
where $T_y(\Delta_K)\not=0$ since $T_0(\Delta_K)=1$ by Lemma 
\ref{lemm:Toddone1}. 
\end{proof}

\begin{rema}
In \cite{HM2} it was shown that a complete non-singular 
simplicial multi-fan satisfying 
\eqref{eq:Tyn+1} is unique up to isomorphisms and is 
isomoprhic to the fan associated to the $n$-dimensional 
projective space $\Proj^n$. There are exactly $n+1$ primitive 
generating vectors $\{v_i\}_{i=1}^{n+1}$ and they satisfy the relation 
\[ v_1+v_2+\cdots +v_{n+1}=0. \]
\end{rema}

\begin{prop}\label{prop:N=n}
Let $\Delta$ be a complete multi-fan of dimension $n$ 
satisfying condition \eqref{eq:Toddone}. If 
$c_1(\Delta)$ is $T$-Cartier divisible by $n$, then 
the following two possibilities occur. 
\begin{enumerate}
\renewcommand{\labelenumi}{(\alph{enumi})} 
\item $\Delta$ is non-singular and 
\[ T_y(\Delta)=(1-y)\sum_{k=0}^{n-1}(-y)^k. \]
\item $n\geq 2$ and $\Delta$ has a unique $K\in \Sigma$ such that 
$\hatH_K\not=\emptyset$. In this case 
\begin{align*}
 T_y(\Delta)&=\sum_{k=0}^{n}(-y)^k, \\
 T_y(\Delta_K)&=\sum_{k=0}^{n-2}(-y)^k.
\end{align*}
and 
\[ \hatT_y(\Delta)=T_y(\Delta)+(-y)T_y(\Delta_K)=
 (1-y)\sum_{k=0}^{n-1}(-y)^k. \]
\end{enumerate}
\end{prop}

\begin{rema}
In \cite{HM2} it was shown that in the case (a) there are 
exactly $n+2$ elements 
in $\sigmone$ and the corresponding primitive generating vectors 
satisfy the relations (under a suitable numbering) 
\[ v_0+v_{n+1}+\sum_{i=2}^nk_iv_i=0,\ v_1+v_2+\cdots +v_n=0,\]
where $(\sum_{i=2}^{n}k_i)+2$ is divisible by $n$. It was also 
shown that $c_1(\Delta)$ of a complete non-singular simplicial 
multi-fan $\Delta$ satisfying such relations is divisible by $n$. 

In the case (b) it will be shown in the proof that there are 
exactly $n+1$ elements in $\sigmone$, and the corresponding 
primitive generating 
vectors satisfy the relations (under a suitable numbering)
\begin{equation*}\label{eq:N=nvi} 
2\sum_{i=1}^{n-1}v_i+v_n+v_{n+1}=0. 
\end{equation*}
\end{rema}

\begin{proof}
We first show that the 
orbifold $T_y$-genus must be of the form
\begin{equation}\label{eq:TyN=n}
 \hatT_y(\Delta)=(1-y)\sum_{k=0}^{n-1}(-y)^k. 
\end{equation}

$\hatT_y(\Delta)$ is a polynomial in $-y$ of 
degree $n$ 
divisible by $\sum_{k=0}^{n-1}(-y)^k$ by Proposition 
\ref{prop:hatTdiv}. We shall show that the constant term 
and the coefficient of the highest 
term are equal to $1$. Then it would prove that 
$\hatT_y(\Delta)$ must be of the form \eqref{eq:TyN=n}. 
The constant term of $T_y(\Delta)$ is 
$h_0(\Delta)=T_0(\Delta)=1$ 
and its coefficient of the highest term is 
$h_n(\Delta)=h_0(\Delta)=1$. So it suffices to show 
that $(-y)^{f_{K,h}}T_y(\Delta_K)$ has no constant term 
and its highest degree is less than $n$ for any $K\in \sigmk$ 
with $k>0$ and $h\in \hatH_K$. 

As $f_{K,h}>0$ by \eqref{eq:fpositive} there is no constant 
term. On the other hand $T_y(\Delta_K)$ is a polynomial of 
degree $n-k$ so that the highest degree of 
$(-y)^{f_{K,h}}T_y(\Delta_K)$ is $f_{K,h}+n-k$ and it is 
less than $n$ by \eqref{eq:fpositive}. This finishes the 
proof of \eqref{eq:TyN=n}.

If $\hatT_y(\Delta)$ is of the form (\ref{eq:TyN=n}), then
$h_{n-1}(\Delta)=1 \ \text{or $2$}$ and 
\begin{equation*}
 e_1(\Delta)=nh_n(\Delta)+h_{n-1}(\Delta)=
 \begin{cases}
 n+1 ,& \quad h_{n-1}(\Delta)=1, \\
 n+2 ,& \quad h_{n-1}(\Delta)=2.
 \end{cases}
\end{equation*}

\begin{claim} Case 1. $h_{n-1}(\Delta)=2,\ e_1(\Delta)=n+2$. We have
\[ e_n(\Delta)=2n,\quad h_k(\Delta)=
 \begin{cases}
 1,\quad &k=0,n, \\
 2,\quad &1\leq k\leq n-1
 \end{cases} \]
\end{claim}
We consider the link $Lk\{i\}$ of $i\in\sigmone$ in $\Sigma$. 
The number $\# Lk\{i\}$ of $Lk\{i\}$ is at least $n$ in general 
by the completeness of $n$-dimensional multi-fan $\Delta$. 
In the present case it is equal to $n$ or $n+1$. 

We show that there is a vertex of $\Sigma$, i.e., an element 
$i\in \sigmone$ with $\# Lk\{i\}=n$. In fact if $\# Lk\{i\}=n+1$ 
for all $i\in \sigmone$, 
then the set $\{i\in \sigmone\}$ would form an $n$-dimensional 
simplex in $\Sigma$. This is a contradiction because the highest 
dimension of simplices in $\Sigma$ is $n-1$, since the dimension 
as a simplex of $I\in \sigmn$ is $n-1$. 

If $\# Lk\{i\}=n$, then the $n$ vertices of $Lk\{i\}$ do not 
form an $(n-1)$-dimensional simplex. For otherwise they together with 
$i$ would form an $n$-dimensional simplex. Hence the star $St\{i\}$ 
of $i$, i.e., the join of $i$ with $Lk\{i\}$ consists of $n$ simplices 
in $\sigmn$. 

Assume now that $\# Lk\{i\}=n$. Then there is a unique $i'\in \sigmone$ 
such that $Lk\{i'\}=Lk\{i\}$. Each of the stars of $i$ and $i'$ consists 
of $n$ simplices in $\sigmn$. Hence 
\[ e_n(\Delta)=2n. \]
We also have 
\begin{equation}\label{eq:ek1} 
 e_k(\Delta)=\begin{cases}
  \binom{n}{k}+2\binom{n}{k-1}, \quad &1\leq k\leq n-1, \\
  \ 2n, \qquad &k=n, \\
  \ \ 1, \qquad  &k=0.
 \end{cases} 
\end{equation} 

On the other hand 
we have the relations 
\[ e_k(\Delta)=\sum_{i=0}^k\binom{n-i}{n-k}h_{n-i}(\Delta). \] 
This follows from
\[ T_y(\Delta)=\sum_{k=0}^nh_k(\Delta)(-y)^k=
 \sum_{k=0}^ne_k(\Delta)(-y-1)^{n-k}. \]
If we put 
\begin{equation}\label{eq:hk} 
h_k(\Delta)=
 \begin{cases}
 1,\quad &k=0,n, \\
 2,\quad &1\leq k\leq n-1
 \end{cases}
\end{equation}
in the above relations, then we have exactly the same values 
\eqref{eq:ek1} for $e_k(\Delta)$. It follows that 
$h_k(\Delta)$ is given by \eqref{eq:hk}. This proves Claim. 

Since $h_k(\Delta)$ is given by \eqref{eq:hk} we have 
\[ T_y(\Delta)=1+2\sum_{k=1}^{n-1}(-y)^k+(-y)^n
 =(1-y)\sum_{k=0}^{n-1}(-y)^k=\hatT_y(\Delta). \] 
We conclude that $\Delta$ is non-singular.

\begin{claim}
Case 2. $h_{n-1}(\Delta)=1,\ e_1(\Delta)=n+1$. We have 
\[ e_n(\Delta)=n+1,\quad h_k(\Delta)=1,\ 0\leq k\leq n. \]
\end{claim}
The fact that $e_1(\Delta)=n+1$ implies that $Lk\{i\}$ 
consists of $n$ vertices for any $i\in \sigmone$. This 
in turn implies that $e_n(\Delta)=n+1$. Then $h_k(\Delta)$ equals 
$1$ for all $0\leq k\leq n$. In fact 
\[ n+1=e_n(\Delta)=\sum_kh_k(\Delta)\geq n+1 \] 
since $h_k(\Delta)>0$. Hence $h_k(\Delta)$ must be equal to $1$ for 
all $k$. This proves Claim. 

Then we see that
\begin{equation}\label{eq:TyN=n2}
 T_y(\Delta)=\sum_{k=0}^nh_k(\Delta)(-y)^k=\sum_{k=0}^n(-y)^k.
\end{equation}
Since $T_y(\Delta)\not= \hatT_y(\Delta)$, 
$\Delta$ can not be non-singular.  

Hereafter we assume that $\sigmone=\{1,2,\ldots,n+1\}$. 
We shall show that there is a unique $K\in \Sigma$ such that 
$\hatH_K\not=\emptyset$. Note that the generating vectors 
$\{v_i\}_{i=1}^{n+1}$ satisfy a relation of the form:
\begin{equation}\label{eq:vn+1}
 \sum_ia_iv_i=0 \ \ \text{with}\ \ a_i\in \Z_{>0}. 
\end{equation}
To see this, we write 
$v_{n+1}$ in the form
\begin{equation*}
 v_{n+1}=\sum_{i=1}^nb_iv_i,\ b_i\in \Q. 
\end{equation*}  
The convexity of rational $n$-dimensional cones 
$C(I),\ I\in \sigmn,$ implies that $b_i$ must be negative 
for all $i$. Hence \eqref{eq:vn+1} follows. We assume that 
the greatest commom divisor of $\{a_i\}$ is equal to $1$. 

We may assume that the lattice $L$ is $\Z^n$. We write 
\[ (v_1,v_2,\ldots,v_n,v_{n+1})=(e_1,e_2,\ldots,e_n)A, \]
with $(n,n+1)$ matrix $A$ where 
$e_i$ is the standard unit vector for $1\leq i\leq n$. 
If $A_i$ denotes the $(n,n)$ matrix obtained from $A$ by 
deleting the $i$-th column, then we see easily that
\[ |\det A_i|=da_i \]
for some positive integer $d$. We shall show later that $d$ is equal 
to $1$. 

We put $\{i\}^*=\{1,2,\ldots,n,n+1\}\setminus \{i\}\in \sigmn$. 
If $L_{\{i\}^*,\V}$ denotes the lattice generated by 
$v_1,\ldots,v_{i-1},v_{i+1},\ldots,v_{n+1}$, then
$H_{\{i\}^*}=L/L_{\{i\}^*,\V}$. Since 
\[ (v_1,\ldots,v_{i-1},v_{i+1},\ldots,v_{n+1})=(e_1,\ldots,e_n)A_i ,\]
we have 
\begin{equation}\label{eq:detAi}
 |H_{\{i\}^*}|=|\det A_i|=da_i. 
\end{equation}

\begin{claim}
\[ d=1,\quad \sum_{i=1}^{n+1}a_i=2n. \]
\end{claim}
\noindent
Proof of Claim. 

\noindent
Step 1. \ We show $d=1$ and 
$\sum_{i=1}^{n+1}a_i\leq 2n$. 

We put $-y=1$ in the polynomial $\hatT_y(\Delta)$ in $-y$. 
By \eqref{eq:TyN=n} we have 
\[ 2n=\hatT_y(\Delta)|_{-y=1}=\sum_{k=0}^n\sum_{K\in \sigmk}
 \sum_{h\in \hatH_K}T_y(\Delta_K)|_{-y=1}. \]
If we put $T_y(\Delta_K)=\sum_{l=0}^{n-k}h_l(\Delta_K)(-y)^l$ for 
$K\in \sigmk$, then Lemma \ref{lemm:Toddone1} and Lemma 
\ref{lemm:hpositive} yield $h_l(\Delta_K)>0$. Therefore 
$T_y(\Delta_K)|_{-y=1}\geq n+1-k$. Hence 
\begin{equation}\label{eq:less2n}
\sum_{k=0}^n\sum_{K\in \sigmk}(n+1-k)|\hatH_K|\leq 2n. 
\end{equation}

On the other hand we have
\begin{equation}\label{eq:less2n2}
 \sum_{k=0}^n\sum_{K\in \sigmk}(n+1-k)|\hatH_K|
 =\sum_{I\in\sigmn}|H_I|. 
\end{equation}
In fact $H_I=\bigsqcup_{K\subset I}\hatH_K$, and the number of 
$I=\{i\}^*\in \sigmn$ containing a fixed $K\in\sigmk$ is equal to 
$n+1-k$. Combining \eqref{eq:less2n2} with \eqref{eq:detAi} 
and \eqref{eq:less2n} 
we obtain 
\begin{equation*}
 d\sum_{i=1}^{n+1}a_i=\sum_{i\in \sigmone}|H_{\{i\}^*}|\leq 2n. 
\end{equation*}
Since all $a_i$ are positive, $d$ must be equal to $1$, and we get
$\sum_{i=1}^{n+1}a_i\leq 2n$. 
\vspace*{0.3cm}\\
Step 2. We show $\sum_{i=1}^{n+1}a_i=2n$. 

Since $c_1(\Delta)$ is $T$-Cartier divisible by $n$, the equivariant
first Chern class is of the form 
\begin{equation}\label{eq:ndivisible}
 c_1^T(\Delta)=\sum_{i\in \sigmone}x_i=nx+u 
\end{equation}
with $u\in H^2(BT)$ and $\iota_I^*(x)\in L^*$ for all $I\in\sigmn$. 
In particular $u^I=\iota^*(c_1^T(\Delta))$ lies in $L^*$. 
Since $\iota_I^*(x)\in L^*$, we have 
\[ \l u^I,v_i\r\equiv\l u,v_i\r \bmod n \]
for all $I\in\sigmn$.

We may suppose that $a_1\geq a_2\geq\cdots\geq a_n\geq a_{n+1}$ without 
loss of generality. Since $\sum_{i=1}^{n+1}a_i=2n$ , we must 
have $a_n=a_{n+1}=1$. 
Put $I=\{n+1\}^*$ and $I'=\{1\}^*$. Then 
\[ -\sum_{i=1}^na_i=\l u^I,v_{n+1}\r \equiv \l u,v_{n+1}\r 
 \equiv \l u^{I'},v_{n+1}\r =1=a_{n+1}\  (\bmod\ n). \]
Therefore $\sum_{i=1}^{n+1}a_i$ is divisible by $n$. 
But $\sum_{i=1}^{n+1}a_i\leq 2n$ and $a_i\geq 1$ for all $i$. 
Hence $\sum_{i=1}^{n+1}a_i=2n$. 
\vspace*{0.3cm}\\
Step 3. We show that $2x_{n+1}$ is $T$-Cartier. 

We first show that
\begin{equation}\label{eq:xn+1}
 c_1^T(\Delta)=2nx_{n+1}+u^I 
\end{equation}
for $I=\{n+1\}^*$.
For that purpose put $x'=c_1^T(\Delta)-u^I$. 
Then $\iota_I^*(x')=0$. Since the kernel of $\iota_I^*$ is 
generated by $x_{n+1}$ we see that $x'=mx_{n+1}$ for 
some $m\in \Z$. Hence $\l\iota_{I'}^*(x'),v_{n+1}\r=m$ for
$I'=\{1\}^*$. On the other hand 
$\l \iota_{I'}^*(c_1^T(\Delta)),v_{n+1}\r=1$, 
and $\l u^I,v_{n+1}\r=-\sum_{i=1}^na_i=-2n+1$. 
Hence we get $m=\l\iota_{I'}^*(x'),v_{n+1}\r=2n$.

Compairing \eqref{eq:xn+1} with \eqref{eq:ndivisible} 
we see that 
\[ \l u-u^I,v_i\r\equiv 0 \bmod n \]
for all $1\leq i\leq n$. Since $|A_{n+1}|=a_{n+1}=1$ 
the collection $\{v_i\}_{i=1}^n$ generates the lattice $L$. It follows that 
$ u-u^I=nu_1 $ for some $u_1\in L^*$. Then 
\[ n(2x_{n+1})=nx+u-u^I=n(x+u_1), \]
and $2x_{v+1}=x+u_1$. Since $x$ is $T$-Cartier by assumption, 
$2x_{v+1}$ is $T$-Cartier. 
\vspace*{0.3cm}\\
Step 4. We shall show that $(a_1,\ldots,a_{n-1},a_n,a_{n+1})
=(2,\ldots,2,1,1)$ 
which will prove \eqref{eq:N=nvi}. 

$v_1,\ldots,v_n$ form a basis of $L$ since $|A_{n+1}|=1$. 
Therefore we may assume that $L$ is $\Z^n$ and $v_i=e_i$, 
the standard unit vector. 
The fact that $2x_{n+1}$ is $T$-Cartier is equivalent to
\[ 2u_{n+1}^I=\iota_I^*(2x_{n+1})\in L^*\ \ \text{for}\ \ 
I=\{i\}^*,\ 1\leq i\leq n. \]
If $\{e_i^*\}$ is the basis of $L^*$ dual to $\{e_i\}$, 
then 
\[ u_{n+1}^I=-\frac{1}{a_i}e_i^* \]
for $I=\{i\}^*$ since $v_{n+1}=-\sum_{i=1}^na_ie_i$. 
Therefore from the above condition $2u_{n+1}^I\in L^*$ 
it follows that $a_i$ must be equal to $1$ or $2$. 
But $\sum_{i=1}^{n+1}a_i=2n$. Hence 
\[ (a_1, \ldots, a_{n-1},a_n,a_{n+1})=(2,\ldots,2,1,1). \]

Now we put $K=\{n,n+1\}\in\Sigma^{(2)}$. 
If $J\in \Sigma$ does not contain $K$, then $J$ is a 
subset of $\{1,\ldots, n-1, n\}$ or $\{1,\ldots, n-1, n+1\}$. 
Since 
$v_1,\ldots, v_{n-1},v_n$ and $v_1,\ldots, v_{n-1},v_{n+1}$ 
are bases of $L$, $H_J$ is a trivial group. 

On the other hand $L_K$ is generated by $v_n$ and $(v_n+v_{n+1})/2$, 
and $L_{K,\{v_n,v_{n+1}\}}$ is generated by $v_n,v_{n+1}$. 
This implies that
$H_K=L_K/L_{K,\{v_n,v_{n+1}\}}\cong\Z/2$. As a representative $v(h)$ 
of its generator $h$ we can take $(v_n+v_{n+1})/2$.
Similarly for $I_i=\{1,\ldots,n,n+1\}\setminus \{i\}$ 
with $1\leq i\leq n-1$, we have 
$H_{I_i}=L/L_{I_i}\cong \Z/2$. This implies that 
$\hatH_K=\{h\}$, and
$K$ is the unique element in $\Sigma$ such that 
$\hatH_K\not=\emptyset$. Moreover if we identify $L^*$ 
with $\Z^n$, then 
\[ u_n^K=-\sum_{i=1}^{n-1}\frac{1}{2(n-1)}e_i+e_n \ \ \text{and}
 \ \ u_{n+1}^K=-\sum_{i=1}^{n-1}\frac{1}{2(n-1)}e_i. \]
It follows that $f_{K,h,i}=\l u_i^K, v(h)\r=\frac12$ for 
$i=n,n+1$, and $f_{K,h}=1$. Therefore 
\[ 
 \hatT_y(\Delta) =T_y(\Delta)+(-y)T_y(\Delta_K). \]
From this and \eqref{eq:TyN=n2} we see also that
\[ T_y(\Delta_K)=\sum_{k=1}^{n-2}(-y)^k. \] 
\end{proof}   

\begin{rema}
Let $\Delta$ be a fan of the form (a) in Proposition 
\ref{prop:N=n}. The corresponding toric variety is a 
projective space bundle over a projective line, c.f. \cite{HM2}
and \cite{Fuj}. In the notation of \cite{Fuj} it is written 
$\Proj_{\Proj^1}(\mathcal{O}(k_1)\oplus\mathcal{O}(k_2)\oplus\cdots 
 \mathcal{O}(k_n))$ with $\sum_{i=1}^nk_i=2$. 
The toric variety corresoponging to the fan of the form (b) 
in Proposition \ref{prop:N=n} is the weighted projective space 
$\Proj^n(2,\ldots ,2,1,1)$. Its orbifold structure is the one 
as the quotient $(\C^{n+1}\setminus \{0\})/\C^*$ where the action 
of $\C^*$ is given by 
\[ z(z_1,\ldots,z_{n-1},z_n,z_{n+1})=
 (z^2z_1,\ldots,z^2z_{n-1},zz_nzz_{n+1}). \]
See \cite{Fuj}. 
It should be noticed that the action ot the finite group 
$(\Z/2)^{n-1}$ on $\Proj^n$ defined by 
\[ (g_1,\ldots,g_{n-1})[z_1,\ldots,z_{n-1},z_n,z_{n+1}]=
 [g_1z_1,\ldots,g_{n-1}z_{n-1},z_n,z_{n+1}],\ g_i=\pm 1, \] 
gives the same algebraic varety $\Proj^n(2,\ldots ,2,1,1)$ but a 
different orbifold structure.
\end{rema}

Since a toric variety is determined by its fan, we 
obtain
\begin{coro}\label{coro:N=n}
Let $X$ be a $\Q$-factorial complete toric variety of dimension $n$ 
and $K_X$ denote the canonical divisor of $X$.  
If there exists a $T$-Cartier divisor $D$ such that $K_X$ is linearly 
equivalent to $nD$, then $X$ is isomorphic to a projective space 
bundle $\Proj_{\Proj^1}(\mathcal{O}(k_1)\oplus\mathcal{O}(k_2)\oplus\cdots 
 \mathcal{O}(k_n))$ with $\sum_{i=1}^nk_i=2$ or to 
$\Proj^n(2,\ldots ,2,1,1)$ as a toric variety. 
\end{coro}

Fujino \cite{Fuj} classified the $n$-dimensional projective toric varieties 
$X$ such that $K_X$ is $\Q$-Cartier and numerically equivalent 
to $nD$ for some Cartier divisor $D$. When $X$ is $\Q$-factorial
the conclusion is the same as Corollary \ref{coro:N=n}. 

\providecommand{\bysame}{\leavevmode\hbox to3em{\hrulefill}\thinspace}

\end{document}